\documentclass[reqno]{amsart}
\usepackage{amsfonts}
\usepackage{amssymb}
\usepackage{amsthm}
\usepackage{upref}
\makeatletter
\def\LaTeX{\leavevmode L\raise.42ex
    \hbox{\kern-.3em\size{\sf@size}{0pt}\selectfont A}\kern-.15em\TeX}
 \makeatother
 
\numberwithin{equation}{section}

\newtheorem{lemma}{Lemma}[section]
\newtheorem{theorem}[lemma]{Theorem}

\newtheorem{proposition}[lemma]{Proposition}

\theoremstyle{definition}

\newtheorem{definition}[lemma]{Definition}

\newtheorem{remark}[lemma]{Remark}

\newcommand{\Det}{\operatorname{Det}}
\newcommand{\tr}{\operatorname{tr}}
\newcommand{\Tr}{\operatorname{Tr}}

\newcommand{\HH}{\mathsf{H}}
 
  \newcommand{\e}{\eqref}

\newcommand{\q}{\quad}

\newcommand{\wt}{\widetilde}

\newcommand{\z}{\zeta}

\renewcommand\Im{\operatorname{Im}}
\renewcommand\Re{\operatorname{Re}}

\newenvironment{pf}{\begin{proof}}{\end{proof}}

\def\qqq{\mathrel{\subset\mkern-15mu\lower.38ex\hbox{${\scriptscriptstyle\rightarrow}$}}}

\let\goth\mathfrak
\let\cal\mathcal

\let\Bbb\mathbb

\begin{document}

\title[trace formula]
{A trace formula   for differential operators \\ of  arbitrary   order}

\author{ J. \"Ostensson}
\address{Department of Mathematics,
Uppsala University, Box 480,
SE-751 06 Uppsala, SWEDEN}
\email{ostensson@math.uu.se}

\author{ D. R. Yafaev}
\address{ IRMAR, Universit\'{e} de Rennes I\\ Campus de
  Beaulieu, 35042 Rennes Cedex, FRANCE}
\email{yafaev@univ-rennes1.fr}
\thanks{The first author is grateful to Ari Laptev for useful and stimulating discussions. The second author was partially supported by the project NONAa, ANR-08-BLANC-0228}


\dedicatory{To the memory of Israel~Cudicovich~Gohberg}

\begin{abstract}
  An operator $H=H_{0}+V$ where $H_{0}=i^{-N} \partial^{N}$ $(N$ is arbitrary) and $V$  is   a differential operator  of order $N-1$ with coefficients decaying   sufficiently rapidly  at infinity    is considered
  in the space $L^2 (\Bbb R)$. The goal of the paper is to find an expression for the trace of the difference of the resolvents $(H-z)^{-1}$ and $(H_{0}-z)^{-1}$ in terms  of the Wronskian of appropriate solutions to the differential equation $H u =z u$.  This also leads to a representation  for the   perturbation determinant of the pair $H_{0}, H$.
\end{abstract}

\maketitle


\section{Introduction}

{\bf 1.1. }
In the framework   of the general operator theory in an abstract Hilbert space, the spectral theory of
differential operators
 \begin{equation}
H=i^{-N} \partial^N +  v_{N }(x) \partial^{N-1}+\cdots +v_{2} (x) \partial +v_1(x),\q \partial=  d /dx,
 \label{eq:H}\end{equation}
 is the same for all values of $N$. However, from the point of view of differential equations the problems are essentially different for $N=2$ (for  $N=1$  it is trivial) and for larger values of $N$.

 Suppose that the coefficients  $v_{j}(x)$, $j=1,\ldots, N$,  decay   sufficiently rapidly  as $|x|\to\infty$, and set $H_{0}=i^{-N}\partial^{N}$. Let
$R_{0}(z)=(H_{0}-z)^{-1}$ and $R(z)=(H-z)^{-1}$ be the  resolvents of the operators $H_{0}$ and $H$ acting in the space $L^2 (\Bbb R)$. The self-adjointness of the operator $H$ is inessential for us, and we do not assume it.

The  main goal of the present paper is to find an expression for the trace
 \begin{equation}
 \Tr\big( R(z)-R_{0} (z)\big)
  \label{eq:trace}\end{equation}
   in terms of solutions to  the differential equation $H u =z u$.
 In the case $N=2$ such an expression   was found     by V. S. Buslaev and L. D. Faddeev in paper  \cite{BF}. They considered   the problem on the half-line, and   the problem on the whole line was discussed by L. D. Faddeev in \cite{Finv}.

 \medskip

{\bf 1.2. }
    Let us introduce the notation
\begin{equation}
 \{u_{1} ,\ldots,   u_{N}  \} =
  \begin{pmatrix}
u_{1}  &\ldots  &u_{N}
 \\
u'_{1}   &\ldots  &u'_{N}
 \\
 \vdots&  \vdots  &  \vdots
  \\
u_{1}^{(N-1)}  &   \ldots  &u^{(N-1)}_{N}
 \end{pmatrix}
 \label{eq:Wrons}\end{equation}
 for the Wronskian matrix of solutions $u_{1} =u_{1}(x,z),\ldots, u_N =u_N  (x,z)$  of the differential equation
\begin{equation}
i^{-N} u^{(N)} (x) +   v_{N} (x)u^{(N-1)}(x) +\cdots +v_2 (x) u'(x)  +v_1(x)u(x)= z u(x).
 \label{eq:Heh}\end{equation}
 We always assume that $z\in{\Bbb C}\setminus[0,\infty)$ if $N$ is even and that $\Im z \neq 0$ if $N$ is odd. Let   $\z_{j}$ be solutions of the equation $\z^{N}=i^N z$. We  suppose that
 \begin{equation}
 \Re \z_{j}>0\q {\rm for}  \q j=1,\ldots, n \q {\rm and} \q \Re \z_{j}<0
 \q {\rm for}  \q j=n+1,\ldots, N.
 \label{eq:zeta}\end{equation}
 Here $n=N/2$ if $N$ is even and $n=(N-1)/2$ for $\Im z>0$ and $n=(N+1)/2$ for $\Im z<0$ if $N$ is odd.

We  first explain our result for the case of functions $v_{j} (x)$ with compact supports. We write $x<\!\!< 0$ if $x$ lies to the left of the supports of all $v_{j} (x)$  and $x>\!\!>0$ if $x$ lies to the right of this set.
    Let $u_{j} (x,z)$ be solutions of
  equation \e{eq:Heh}  such that
  \begin{equation}
u_{j} (x,z)=e^{\z_{j} x }\; {\rm for}\; x<\!\!<0 \; {\rm if}\; j=1,\ldots,n \q  {\rm and}\q
{\rm for}\; x>\!\!>0 \; {\rm if}\; j=n+1,\ldots,N.
 \label{eq:solu}\end{equation}
   Let
 \begin{equation}
 {\sf W} (x, z)= \det\{u_{1} (x,z),\ldots, u_{n} (x,z), u_{n+1} (x,z), \ldots, u_{N} (x,z)\}
 \label{eq:WRO}\end{equation}
 be  the determinant of matrix \e{eq:Wrons},  and let
 \begin{equation}
 {\sf W}_{0} (z)= \det\{e^{\z_{1} x},\ldots, e^{\z_{n} x}, e^{\z_{n+1} x}, \ldots, e^{\z_{N} x} \}
 \label{eq:WRO1}\end{equation}
be  the corresponding Wronskian  for the  ``free" case where $v_{j}=0$ for all $j=1,\ldots, N$. Of course,
\begin{equation}
 {\sf W} (x_{2}, z)= \exp\big( -i^N \int_{x_1}^{x_{2}} v_{N}(y) dy\big) {\sf W} (x_1, z)
 \label{eq:WROx}\end{equation}
 for arbitrary points $x_{1}$ and $x_2$.
We emphasize that the    Wronskians ${\sf W} (x, z)$ and $ {\sf W}_{0} (z)$ depend on the order of numeration of the numbers $\z_{j}$, but  the normalized  Wronskian
 \begin{equation}
\Delta(x, z)={\sf W}  (x, z)/ {\sf W} _{0}(z)
\label{eq:det2}\end{equation} does not depend on it.

Our main result is that the normalized  Wronskian
 satisfies (for all $x$ and all regular points $z$ of the operator $H$) the equation
\begin{equation}
\Tr\big( R(z)-R_{0}(z)\big)= - \Delta(x, z)^{-1} d\Delta(x, z)/dz,
\label{eq:tr}\end{equation}
which we call the trace formula in this paper.  Thus the trace of the difference of the resolvents admits an explicit expression in terms of properly chosen solutions of equation \e{eq:Heh}.

Then we extend  representation \e{eq:tr} to general short-range coefficients $v_{j}(x)$ satisfying the assumption
 \begin{equation}
\int_{-\infty}^\infty | v_{j}(x)|^2 (1+ x^2)^\alpha dx<\infty,\q \alpha>1/2,\q j=1,\ldots, N,
\label{eq:V}\end{equation}
 only. In this case the functions $u_{j}(x,z)$ in definition
\e{eq:WRO} are the solutions of equation \e{eq:Heh}  such that
\begin{equation}
u_{j} (x,z)=e^{\z_{j} x }(1+o(1))
\label{eq:as10}\end{equation}
    as $x\to-\infty$ if $j=1,\ldots,n$ and   as $x\to+ \infty$ if $j=n+1,\ldots, N$. Here and below {\it all asymptotic relations for solutions of equation \e{eq:Heh} are supposed to be $N-1$ times differentiable in } $x$.
    We emphasize that for $N > 2$  asymptotics \e{eq:as10} DO NOT determine the solutions   of equation \e{eq:Heh}  uniquely. However, the Wronskian \e{eq:WRO}  does not depend on specific choice of the solutions satisfying \e{eq:as10}. Thus we do not need the construction of the book \cite{BDT}  by R. Beals, P. Deift and C. Tomei devoted to the inverse scattering problem. In \cite{BDT}  solutions   of equation \e{eq:Heh} were distinguished uniquely (away from some exceptional set of values of $z$) by conditions at both infinities.
        Our construction of solutions of equation \e{eq:Heh} with asymptotics \e{eq:as10} relies on integral equations which are Volterra   equations for $N=2$ but are only Fredholm equations in the general case. Nevertheless for the construction of solutions   with asymptotics \e{eq:as10} as $x\to +\infty$ (as   $x\to -\infty$) we impose conditions on the coefficients $v_{j}(x)$ also as $x\to +\infty$ (as   $x\to -\infty$) only.

Suppose that $v_{N}=0$. Then $ {\sf W} (x,z)= {\sf W} (  z)$ and hence $\Delta (x,z)= \Delta (  z)$ do not depend on $x$.
In this case we identify $\Delta(z)$ with   the perturbation determinant for the pair of  operators $H_{0} $,  $H$. We refer to the book \cite{GoKr} by I. C. Gohberg and M. G. Kre\u{\i}n for a comprehensive discussion of different properties of  perturbation determinants.   Set
 \begin{equation}
V=H-H_{0}=  v_{N } (x) \partial^{N-1}+\cdots +v_{2} (x) \partial +v_{1}(x).
\label{eq:V1}\end{equation}
If $v_{N}=0$, then the operator $V R_{0}(z)$ for $\Im  z \neq 0$  belongs to the trace class ${\goth S}_{1}$, and hence the perturbation determinant
\begin{equation}
D(z)=\Det\big(I+VR_{0}(z)\big)
\label{eq:PD}\end{equation}
is well defined.    Of particular importance is the abstract trace formula
\begin{equation}
\Tr \big(R(z)-R_{0}(z)\big)= -D(z) ^{-1}d D(z) /dz ,
\label{eq:PD1}\end{equation}
which for definition \e{eq:PD} is a direct consequence of the formula for the derivative of a determinant.
Comparing equations \e{eq:tr} and \e{eq:PD1} and using that $\Delta (z)  \to 1$ as $| \Im z| \to \infty$, we   show that
\begin{equation}
 \Det\big(I+V R_{0}(z)\big) = \Delta(z).
\label{eq:PDd}\end{equation}
Thus the perturbation determinant admits an explicit expression in terms of  solutions of equation \e{eq:Heh}.

  If    $v_{N}\neq 0$, then  under assumption  \e{eq:V}  it is still true
 that (for all regular points $z$)
  \begin{equation}
 R(z)-R_{0}(z) \in{\goth S}_{1},
   \label{eq:TSS}\end{equation}
  although   $VR_{0}(z)\not\in{\goth S}_{1}$.
     Without the condition   $v_{N }= 0$,   equation  \e{eq:PD1} is satisfied for   so called  generalized perturbation determinants $\wt{D}(z)$ which are defined up to constant factors (see subs.~6.2). According to equation  \e{eq:tr} in the general case for every fixed $x\in \Bbb R$,  the function  $\Delta(x, z)$ differs from each generalized perturbation determinant by a  constant (not depending on $z$) factor.

   \medskip

{\bf 1.3.}
A preliminary step in the proof of the trace formula  \e{eq:tr} is to find a convenient representation for the resolvent $R(z) $ of the operator $H$.
This construction   goes   probably back to the beginning of the twentieth century. We refer to relatively recent books
 \cite{AhGl,BDT,Nai} where its different versions can be found. We start, however, with writing down necessary formulas in  a form convenient for us.

  A differential equation of order $N$ can, of course,  be rewritten as a special system of $N$ differential equations of the first order. A consideration of first order systems without special assumptions on their coefficients gives more general and transparent results.
 A large part of the paper is written in terms of solutions of first order systems which implies the results about solutions of differential equations of an arbitrary order as their special cases.

Let us briefly discuss the structure of the paper. In Sections~2 and 3 we collect necessary formulas for
 solutions of first order systems. They are used in Section~4 for the construction of the integral  kernel
 $R(x,y, z) $ of   $R(z) $. In particular, we obtain a new representation for the integral
 \begin{equation}
 \int_{x_1}^{x_{2}} R(y,y, z)dy
 \label{eq:Di}\end{equation}
 where the points $x_{1}, x_{2}\in{\Bbb R}$ are arbitrary.
 Then passing to the limit $x_{1}\to-\infty$, $x_{2}\to+\infty$, we   prove the trace formula  \e{eq:tr} for   the coefficients $v_{j}$, $j= 1,\ldots,  N$, with compact supports.
   A construction of solutions of equation   \e{eq:Heh}  with asymptotics \e{eq:as10}  is given in Section~5. Here we again   first consider a general system of $N$ differential equations of the first order. Finally, in Section~6 we give the definition of the normalized Wronskian for  operators $H$ with arbitrary short-range coefficients and extend the trace formula  to the general case. At the end we prove that the normalized Wronskian coincides with the perturbation determinant.

  \medskip

{\bf 1.4.}
We note that there exists a somewhat different approach to proofs of formulas of type
\e{eq:PDd}. It consists of a direct calculation of determinant \e{eq:PD} whereas we proceed from a calculation of trace \e{eq:trace}. In this way formula \e{eq:PDd} was proven in \cite{JP} for the Schr\"odinger operator on the half-line. In  \cite{JP} the Fredholm expansion of determinants was used.

A general approach to a calculation of determinants $\Det\big(I+K\big)$ was proposed in the book  \cite{GGK} by I. C. Gohberg, S. Goldberg and N. Krupnik. In this book integral operators $K$ with so called semi-separable kernels were considered. It is important that operators $K=V R_{0}(z)$ fit into this class. This approach was applied to the Schr\"odinger operator in paper \cite{GeMa}.

The authors thank F. Gesztesy for pointing out references \cite{JP, GGK, GeMa}.

\section{Resolvent kernel}

In this section we consider an auxiliary vector problem.

 \medskip

{\bf 2.1. }
Suppose that the eigenvalues $\z_{j}$, $j=1,\dots, N$,  of an $N\times N$ matrix ${\bf L}_{0}$ are distinct.  We denote by ${\bf p}_{j}=(p_{1,j},p_{2,j},\ldots, p_{N, j})^t$ (this notation means that the vector  ${\bf p}_{j}$ is considered as a column)  eigenvectors of  ${\bf L}_{0}$  corresponding to its eigenvalues $\z_{j}$ and by ${\bf p}_{j}^*$   eigenvectors of  ${\bf L}_{0}^*$  corresponding to its eigenvalues $\bar{\z}_{j}$.
Recall that $\langle {\bf p}_{j},{\bf p}_{k}^*\rangle=0$ if $j\neq k$ (here $\langle \cdot, \cdot\rangle$ is the scalar product in ${\Bbb C}^N$).
Normalizations of ${\bf p}_{j}$ and ${\bf p}_{j}^*$ are inessential, but we suppose that
$\langle {\bf p}_{j},{\bf p}_{j}^*\rangle =1$. Then the bases  ${\bf p}_{j}$ and ${\bf p}_{j}^*$, $j=1,\dots, N$,  are dual to each other.

 Assume that an $N\times N$ matrix ${\bf V} (x)$ where $x\in {\Bbb R}$ belongs locally to $L^1$ and has compact support. We write $x<\!\!< 0$ if $x$ lies to the left of the support of  ${\bf V} (x)$  and $x>\!\!>0$ if $x$ lies to the right of this set. We put
\begin{equation}
{\bf L} (x)={\bf L}_{0}+{\bf V} (x).
\label{eq:LL}\end{equation}
Consider the homogeneous equation
\begin{equation}
{\bf u}' (x)= {\bf L}(x)  {\bf u} (x)
\label{eq:H2}\end{equation}
for the vector-valued function ${\bf u}(x) =(u_{1 }(x),\ldots, u_N (x))^t$.
For    arbitrary  linearly independent solutions ${\bf u}_{j}(x) =(u_{1,j}(x) ,\ldots, u_{N,j}(x) )^t$ of this equation,  we denote by
\begin{equation}
{\bf U}  (x)=  \begin{pmatrix}
u_{1,1}(x) & u_{1,2} (x)& \ldots&  u_{1,N}(x)
 \\
 u_{2,1}(x) & u_{2,2}(x) & \ldots&  u_{2,N}(x)
 \\
 \vdots&  \vdots& \ddots&\vdots
  \\
u_{N,1}(x) & u_{N,2}(x) & \ldots&  u_{N,N} (x)
 \end{pmatrix} =:\{{\bf u}_{1}(x),{\bf u}_{2}(x),\ldots, {\bf u}_{N}(x)\}
 \label{eq:U}\end{equation}
  the corresponding fundamental matrix. It satisfies the matrix equation
\begin{equation}
  {\bf U}'  (x) = {\bf L}(x)   {\bf U}  (x) .
\label{eq:H3}\end{equation}
It follows that
\begin{align}
d \det {\bf U}  (x)/ d x=&
\det {\bf U}  (x )  \tr \big(  {\bf U}'  (x )  {\bf U}^{-1}(x )\big)
\nonumber\\
= &\det {\bf U}  (x ) \tr {\bf L}(x )
\label{eq:Hdet}\end{align}
and hence
\begin{equation}
\det {\bf U} (x_{2})= \exp\big(   \int_{x_1}^{x_{2}} \tr {\bf L}(y) dy\big) \det {\bf U} (x_1)
 \label{eq:WROxm}\end{equation}
 for arbitrary points $x_{1}$ and $x_2$. Of course $\det {\bf U} (x )\neq 0$ for   all  $x\in \Bbb R$.

 We always suppose that   $\kappa_{j}:=\Re\z_{j}\neq 0$ for all $j=1,\dots, N$.
 Let $n$ and $N-n$ be the numbers of  eigenvalues $\z_{j}$ of the matrix ${\bf L}_{0}$ lying in the right and left half-planes, respectively.
 The cases $n=0$ or $n=N$   where all $\z_{j}$ lie in one of the   half-planes are not excluded.
  Let   ${\bf u}_{j} (x) $ be   solutions of equation \e{eq:H2}  distinguished by the condition
  \begin{equation}
{\bf u}_{j} (x)=e^{\z_{j} x }{\bf p}_{j}\; {\rm for}\; x<\!\!<0 \; {\rm if}\;  \kappa_{j} >0\q  {\rm and}\q
{\rm for}\; x>\!\!>0 \; {\rm if}\;  \kappa_{j}<0 .
 \label{eq:Solu}\end{equation}
     We  denote by ${\sf K}_{+}$ and ${\sf K}_{-}$ the linear spaces  spanned by all solutions ${\bf u}_{j} (x) $ such that  $\kappa_{j}>0 $ and such that  $\kappa_{j}<0 $, respectively. Clearly, $\dim {\sf K}_{+}=n$  and  $\dim {\sf K}_{-}=N-n$.
     We    assume that
 \begin{equation}
{\sf K}_{+}\cap {\sf K}_{-}=\{0\}.
\label{eq:zero}\end{equation}
Then  all nontrivial solutions of equation \e{eq:H2} exponentially grow either as $x\to +\infty$ or  as $x\to -\infty$. In particular, equation \e{eq:H2}
does not have nontrivial solutions ${\bf u} \in L^2 ({\Bbb R} ; {\Bbb C}^N )$.

If ${\bf u}_1 (x) ,\ldots, {\bf u}_n (x)$ and ${\bf u}_{n+1}(x) ,\ldots, {\bf u}_N (x)$ are  {\it  arbitrary}  linear independent solutions
from ${\sf K}_{+}$ and ${\sf K}_{-}$  respectively, then in view of condition \e{eq:zero} all these solutions   are linearly independent. It is now convenient to accept the following

\begin{definition}\label{adm}
 Suppose that $n$ columns of matrix \e{eq:U} form a basis in the linear space ${\sf K}_{+}$ and other $N-n$ columns   form a basis in   ${\sf K}_{-}$. Then the fundamental matrix ${\bf U}  (x) $  is called {\it admissible}.
\end{definition}

Observe that for the ``free" case where   ${\bf V} (x)=0$, we can set
 \begin{equation}
{\bf U}_{0} (x)=\{{\bf p}_{1} e^{\z_{1} x},\ldots, {\bf p}_{N} e^{\z_N x}\}
\label{eq:Hdc}\end{equation}
and
\begin{equation}
 {\sf W}_{0} (x)= \det {\bf U}_{0}(x)= \det \{{\bf p}_{1}   ,\ldots,   {\bf p}_{N} \} \exp\big(\tr {\bf L}_{0}x \big)
\label{eq:WW}\end{equation}
because $\tr {\bf L}_{0}=\sum_{j=1}^N \z_{j}$. Note that $ \det \{{\bf p}_{1}   ,\ldots,   {\bf p}_{N} \} \neq 0$ since all eigenvalues of the matrix $ {\bf L}_{0}$ are distinct. The inverse matrix ${\bf G}_{0}  (x)={\bf U}_{0}^{-1} (x)$ satisfies the relation
\begin{equation}
{\bf G}_{0}^* (x)=\{{\bf p}_{1}^* e^{-\bar{\z}_{1}x} ,\ldots, {\bf p}_{N}^* e^{-\bar{\z}_{N}x}\}.\label{eq:GDd}\end{equation}

  \medskip

 {\bf 2.2. }
 Next we consider the nonhomogeneous equation
\begin{equation}
{\pmb\varphi}' (x)= {\bf L}(x) {\pmb\varphi} (x)  +    {\bf f} (x),\q {\bf f}(x) =(f_{1 }(x),\ldots, f_N (x))^t,
\label{eq:H1}\end{equation}
where the vector-valued  function ${\bf f} (x)$ has compact support.
Let us use the standard method of variation of arbitrary constants
and set
\[
{\pmb\varphi} (x)= {\bf U}  (x)  {\bf q} (x) , \q {\bf q} (x)= (q_{1} (x),\ldots,q_{N} (x))^t,
\]
so that
\begin{equation}
{\pmb\varphi} (x)= \sum_{j=1}^N q_{j} (x) {\bf u}_{j}  (x)  .
\label{eq:H4B}\end{equation}
 Here  ${\bf U}  (x)$   is an arbitrary   admissible fundamental matrix \e{eq:U}.
Then it follows from equation \e{eq:H3} that
\begin{equation}
{\bf q}' (x)= {\bf g} (x) \q {\rm where} \q {\bf g} (x)= {\bf G}(x)  {\bf f} (x)\; {\rm and} \;  {\bf G}(x)={\bf U} ^{-1}(x).
\label{eq:H5}\end{equation}

We are looking for a solution of equation \e{eq:H1}  decaying (exponentially)  as $|x|\to\infty$.  It is convenient to accept convention \e{eq:zeta} on the  eigenvalues $\z_{j} $ of the matrix ${\bf L}_{0}$.  Set
\begin{equation}
\rho_{+}=\min_{j=1,\ldots,n} \Re\z_{j}, \q   \rho_{-}=\min_{j=n+1,\ldots,N} |\Re\z_{j}|
\label{eq:min}\end{equation}
and observe that   estimates
\[
{\bf u}_{j}(x ) = O( e^{-\rho_{\pm}|x|} ),\q x \to \mp\infty,
\]
hold for $j=1, \ldots, n$   and the upper sign as well as for     $j=n+1, \ldots, N$ and the lower sign. Taking into account \e{eq:H4B}, we see that we have to solve equation \e{eq:H5} for different components $q_{j}(x)$ of ${\bf q} (x)$ by different formulas. Namely, we set
\[
\begin{split}
{ q}_{j} (x)&= - \int_{x}^\infty g_{j}(y)    dy,  \q j=1,\ldots, n,
\\
{  q}_{j} (x)&=  \int^{x}_{-\infty }g_{j}(y)    dy,  \q j=n+1,\ldots, N,
\end{split}
\]
where $g_{j} (x)$ are components of $ {\bf g} (x) $.
This leads to the following result.

\begin{proposition}\label{res}
Let assumption  \e{eq:zero} hold, and let  \e{eq:U} be an arbitrary admissible fundamental matrix. Then the function
\begin{equation}
{\pmb\varphi} (x)= - \sum_{j=1}^n {\bf u}_{j}(x)  \int_{x}^\infty  ({\bf G}(y)  {\bf f} (y))_{j}dy
+  \sum_{j=n+1}^N {\bf u}_{j}(x)  \int_{-\infty}^x  ({\bf G}(y)  {\bf f} (y))_{j}dy
\label{eq:res}\end{equation}
satisfies equation \e{eq:H1} and
${\pmb\varphi} (x)=O(e^{-\rho_{\pm} |x|})$   as $x \to\mp \infty$.
\end{proposition}

 Formula \e{eq:res} can be rewritten as
\begin{equation}
{\pmb\varphi} (x)=  \int_{-\infty}^\infty {\bf R}(x,y) {\bf f}  (y) dy
\label{eq:resX}\end{equation}
where the matrix-valued resolvent kernel (or the Green function) ${\bf R}(x,y)=\{ R_{k,l}(x,y)\}$ is defined by the equality
\begin{equation}
  R_{k,l}(x,y)= -\sum_{j=1}^n u_{k,j}(x) g_{j,l}(y)  \theta (y-x)
  + \sum_{j=n+1}^N u_{k,j}(x) g_{j,l}(y)  \theta (x-y).
\label{eq:resZ}\end{equation}
Here $\theta$ is the Heaviside function, i.e., $\theta(x)=1$ for $x\geq 0$ and $\theta(x)=0$ for $x< 0$, and $g_{j,l}$ are elements of the matrix  ${\bf G}$. In the matrix notation  formula \e{eq:resZ} means that
\begin{equation}
{\bf R}(x,y)= -{\bf U} (x) {\sf   P}_{+} {\bf U}^{-1} (y) \theta (y-x)
+   {\bf U} (x) {\sf   P}_{-} {\bf U}^{-1} (y) \theta (x-y),
\label{eq:res1}\end{equation}
where    the  projections ${\sf   P}_\pm $ are defined in the representation ${\Bbb C}^N ={\Bbb C}^n  \oplus {\Bbb C}^{N-n}$ by the block matrices
\[
{\sf P}_+=
\begin{pmatrix}
I_{n} & 0
 \\
0 & 0
 \end{pmatrix},\q
 {\sf    P}_-=
 \begin{pmatrix}
0 & 0
 \\
0 & I_{N-n}
 \end{pmatrix}  .
\]

  Expressions \e{eq:resZ} or \e{eq:res1} do  not of course depend on the choice of   bases in the spaces  $ {\sf K}_{+}$ and $ {\sf K}_{-}$. Indeed, if we choose other bases $\breve{\bf u}_{1}(x),\ldots, \breve{\bf u}_{n}(x)$ and $\breve{\bf u}_{n+1}(x),\ldots, \breve{\bf u}_{N}(x)$, then the corresponding admissible fundamental matrices ${\bf U} (x)$ and $\breve{\bf U} (x)$ are related by the formula   $\breve{\bf U} (x)={\bf U} (x) {\sf F}$ where the operator ${\sf F}: {\Bbb C}^N \to {\Bbb C}^N $ commutes with the projections ${\sf  P}_{\pm}$. It follows that $\breve{\bf U} (x) {\sf  P}_{\pm}\breve{\bf U}^{-1} (y)={\bf U} (x) {\sf  P}_{\pm} {\bf U}^{-1} (y)$.

 \medskip

  Evidently,   the resolvent kernel  \e{eq:res1} is a continuous function of $x$ and $y$ away from the diagonal $x=y$ and
\[
\begin{split}
{\bf R} (x,x+0, z)  &= -{\bf U} (x) {\sf P}_{+} {\bf U}^{-1} (x),
\\
{\bf R} (x,x-0, z)  &= {\bf U} (x) {\sf P}_{-} {\bf U}^{-1} (x).
\end{split}
\]
It follows that
\begin{equation}
{\bf R} (x,x- 0, z)  -
{\bf R} (x,x+ 0,  z)= {\bf I},
\label{eq:Hdc1}\end{equation}
  where ${\bf I}$ is the $N\times N$ identity matrix.

 \medskip

 {\bf 2.3. }
 The results of the previous subsection admit a simple operator interpretation. Consider the space $L^2 ({\Bbb R}; {\Bbb C}^N)$ and define the operator ${\bf H}_{0}$ on the Sobolev class ${\HH}^{1} ({\Bbb R}; {\Bbb C}^N)$  by the formula
 \[
 {\bf H}_{0}=\partial {\bf I} - {\bf L}_{0},\q \partial= d/dx.
 \]
 If ${\bf u}(x)=  u(x) {\bf p}_{j}    $ where $u\in {\HH}^{1} ({\Bbb R} )$, then
 $( {\bf H}_{0}{\bf u})(x)=   (u'(x)-\z_{j } u(x) ) {\bf p}_j$,  and hence the operator $ {\bf H}_{0}$ is {\it linearly} equivalent to a direct sum of the operators of multiplication by $i\xi-\z_{j}$, $\xi\in {\Bbb R}$, $j=1,\ldots, N$, acting in the space $L^2 ({\Bbb R} )$. It follows that the spectrum of the operator $ {\bf H}_{0}$ consists of straight lines passing through all points $-\z_{j}$ and parallel to the imaginary axis.  In particular,  the inverse operator $ {\bf H}_{0} ^{-1}$  exists and is bounded.

 To define the operator
 \[
 {\bf H} =\partial {\bf I} - {\bf L}_{0}-{\bf V} (x),
 \]
 we need the following well known assertion  (see paper \cite {Bia} by M.~Sh.~Birman).

   \begin{lemma}\label{RC}
    Let   $T: L^2 ({\Bbb R}; dx)\to L^2 ({\Bbb R}; d\xi)$ be an integral operator with kernel
\begin{equation}
t(\xi,x)= b(\xi) e^{- i x\xi}v(x).
 \label{eq:TrCl2}\end{equation}
 If $b(\xi)= (\xi^2+1)^{-1/2}$ and
   \begin{equation}
 \lim_{|x|\to\infty}\int_{x}^{x+1} |v (y)|^2 dy=0,
\label{eq:rc}\end{equation}
  then the operator $T$ is compact.
\end{lemma}

If the coefficients of the matrix ${\bf V}(x)$ satisfy condition \e{eq:rc}, then
according to Lemma~\ref{RC} the operator ${\bf V}{\bf H}_{0}^{-1}$ is compact. Hence the operator ${\bf H}$ is closed on ${\HH}^{1} ({\Bbb R}; {\Bbb C}^N)$ and by virtue of the Weyl theorem essential spectra of the operators  $ {\bf H}$ and $ {\bf H}_{0}$ coincide.
 Condition  \e{eq:zero}  implies  that  $0$ is not an eigenvalue of $ {\bf H}$ so that  the inverse operator $ {\bf H} ^{-1}$  exists and is bounded. If the matrix-valued function ${\bf V}(x)$  has compact support, then according to Proposition~\ref{res} the integral kernel of the operator $ {\bf H} ^{-1}$ is given by formula  \e{eq:res1}.

\medskip

{\bf 2.4. }
Let the solutions $ {\bf u}_{j} (x)$
 of equation  \e{eq:H2}  be distinguished by conditions  \e{eq:Solu}.
Let us give expressions for the Wronskian ${\sf W} (x): = \det {\bf U}(x)$ in terms of
 transition matrices ${\bf T}_{\pm}  $ defined as follows.
For $j=1,\ldots, n$    and $x>\!\!> 0$ or  $j=n+1,\ldots, N$    and $x <\!\!< 0$,  we have
\begin{equation}
 {\bf u}_{j} (x) = \sum_{k=1}^{N} t_{j,k}  {\bf p}_{k} e^{\z_kx}
\label{eq:AS}\end{equation}
with some coefficients $t_{j,k} $. Set
\begin{equation}
{\bf T}_{+}   =  \begin{pmatrix}
t_{1,1} & t_{1,2} & \ldots&  t_{1,n}
 \\
t_{2,1} & t_{2,2} & \ldots&  t_{2,n}
 \\
 \vdots&  \vdots& \ddots&\vdots
  \\
t_{n,1} & t_{n,2} & \ldots&  t_{n,n}
 \end{pmatrix}
\label{eq:U2}\end{equation}
and
\begin{equation}
 {\bf T}_{-}   =  \begin{pmatrix}
t_{n+1,n+1}& t_{n+1,n+2} & \ldots&  t_{n+1,N}
 \\
t_{n+2,n+1} & t_{n+2,n+2} & \ldots&  t_{n+2, N}
 \\
 \vdots&  \vdots& \ddots&\vdots
  \\
t_{N,n+1} & t_{N,n+2} & \ldots&  t_{N,N}
 \end{pmatrix}.
\label{eq:U2-}\end{equation}

Consider, for example, ${\bf T}_{+}  $.  Using expressions \e{eq:AS} for $j=1,\ldots, n$, we see that
for $x>\!\!> 0$    matrix \e{eq:U}  equals
 \begin{equation}
{\bf U}(x) =\big\{\sum_{k=1}^{N} t_{1,k} {\bf p}_{k} e^{\z_kx}, \ldots, \sum_{k=1}^{N} t_{n,k}  {\bf p}_{k} e^{\z_kx},
{\bf p}_{n+1}  e^{\z_{n+1} x}, \ldots, {\bf p}_{N}  e^{\z_N x}\big\} .
\label{eq:as2}\end{equation}
Below,  by the calculation of determinants of matrices,  we systematically use that one can add to each column another column multiplied by an arbitrary number. In particular, we have
 \begin{align}
 {\sf W} (x) =&\det\big\{\sum_{k =1}^{n} t_{1,k} {\bf p}_{k}  e^{\z_kx}, \ldots, \sum_{k=1}^{n} t_{n,k} {\bf p}_{k}  e^{\z_kx},
{\bf p}_{n+1} e^{\z_{n+1} x}, \ldots, {\bf p}_{N}  e^{\z_N x}\big\}
\nonumber\\
=& \det  {\bf T}_{+} \,  {\sf W}_{0} (x),\q x>\!\!> 0,
\label{eq:U1}\end{align}
where the free Wronskian $ {\sf W}_{0} (x) $ is given by formula \e{eq:WW}.
In view of relation \e{eq:WROxm}, it follows that for all $x\in{\Bbb R}$
\begin{equation}
 {\sf W} (x) =\exp\big( \tr{\bf L}_{0}x - \int_{x}^\infty \tr {\bf V}(y) dy \big) \det  {\bf T}_{+} \det \{{\bf p}_{1}   ,\ldots,   {\bf p}_{N} \} .
\label{eq:DET}\end{equation}

Quite similarly,   using expressions \e{eq:AS} for $j=n+1,\ldots, N$ and $x <\!\!< 0$, we obtain that
\begin{equation}
 {\sf W} (x) = \det  {\bf T}_{-} \,  {\sf W}_{0} (x),\q x<\!\! <0,
\label{eq:U1-}\end{equation}
and
\begin{equation}
 {\sf W} (x) =\exp\big( \tr{\bf L}_{0}x  + \int^{x}_{-\infty}  \tr {\bf V}(y) dy \big) \det  {\bf T}_{-} \det \{{\bf p}_{1}   ,\ldots,   {\bf p}_{N} \}, \q \forall x\in{\Bbb R}.
\label{eq:DET-}\end{equation}

This leads to the following result.

\begin{proposition}\label{WRTR}
Let ${\bf U}(x)$ be the fundamental matrix \e{eq:U} where   $ {\bf u}_{j} (x)$ are the solutions of equation  \e{eq:H2}  satisfying conditions  \e{eq:Solu}.
Let the transition matrices ${\bf T}_\pm  $ be defined by  formulas \e{eq:AS} -- \e{eq:U2-}. Then the   Wronskian ${\sf W} (x)  = \det {\bf U}(x)$  admits representations \e{eq:DET} and \e{eq:DET-}.
\end{proposition}

Putting together equalities \e{eq:DET} and \e{eq:DET-}, we see that
\[
\det  {\bf T}_+= \exp\big(\! \int_{-\infty}^\infty  \tr {\bf V}(y)dy \big)\det  {\bf T}_{-} .
\]
Assumptions \e{eq:zero}  and $\det  {\bf T}_\pm\neq 0$ are  of course equivalent.

\section{Dual problem}

Some properties of admissible fundamental matrices become more transparent if one
 considers the dual problem corresponding to the matrix-valued function
 \[
 \widetilde{\bf L}(x)=-{\bf L}^* (x).
 \]

 \medskip

 {\bf 3.1. }
  It follows from   equation \e{eq:H3} for $ {\bf U}  (x )$ that the inverse operator
  $  {\bf G}  (x )={\bf U}^{-1}  (x ) $
   satisfies the equation
\begin{equation}
{\bf G}'  (x) = -{\bf G}  (x ) {\bf U} ' (x ) {\bf G}  (x ) =
-{\bf G}  (x ){\bf L}(x )
\label{eq:HT}\end{equation}
which yields the equation
\begin{equation}
\widetilde{\bf U}'  (x) =\widetilde{\bf L} (x)\widetilde{\bf U}  (x)
\label{eq:HTd}\end{equation}
  for the matrix-valued function $\widetilde{\bf U}(x): = {\bf U}^* (x)^{-1}$.
Clearly, $\det\widetilde{\bf U}  (x)\neq 0$ so that $\widetilde{\bf U}(x)$ is a fundamental matrix for this equation. Proposition~\ref{dual} below shows that it is admissible. Set $ \widetilde{\bf U} (x)=\{\tilde{u}_{j,l}(x)\}$, $ {\bf G} (x)=\{g_{j,l}(x)\}$.    We use below that
 \begin{equation}
\overline{\tilde{u}_{l,j}(x)}=g_{j, l}(x) =(-1)^{j+l}m_{ l,j}(x )/ {\sf W}(x ),
\label{eq:m2}\end{equation}
where   $m_{ l,j}(x )$ is the minor of the matrix ${\bf U}(x )$ which is the determinant of the matrix cut down from ${\bf U}(x )$ by removing the row with index $l$ and the column with index $j$.

\begin{proposition}\label{dual}
Let ${\bf u}_{j}(x)  $ be  arbitrary  linear independent solutions of equation  \e{eq:H2}
from ${\sf K}_{+}$ for $j=1, \ldots, n$ and from ${\sf K}_{-}$   for $j=n+1, \ldots, N$, and let $ {\bf U}(x)$ be the corresponding admissible fundamental matrix \e{eq:U}.
Then for all $l=1,\ldots, N$ we have
\begin{equation}
\begin{split}
g_{j,l}  (x)  &= O(e^{ \rho_{-} x}),\q x\to -\infty,\q j=n+1, \ldots, N,
\\
g_{j,l}  (x)  &= O(e^{- \rho_{+} x}),\q x\to +\infty,\q j= 1, \ldots, n ,
\end{split}
 \label{eq:as3c}\end{equation}
with  positive numbers $\rho_{\pm} $   defined in \e{eq:min}.
\end{proposition}

\begin{pf}
Let us prove, for example, the first of relations  \e{eq:as3c}.
  Changing if necessary the numeration, we can suppose that  $j=N$ which is notationally convenient.
For $x<\!\!< 0$ and some numbers $ c_{j,k}$, we have
\[
{\bf U}(x)= \big\{\sum_{k =1}^{n} c_{1,k}  {\bf p}  _{k}e^{\z_kx}, \ldots,
\sum_{k=1}^{n} c_{n,k}  {\bf p} _{k} e^{\z_kx},
\sum_{k =1}^{N} c_{n+1,k}  {\bf p} _{k}e^{\z_kx},\dots, \sum_{k =1}^{N} c_{N,k}  {\bf p} _{k}e^{\z_kx}\big\}.
\]
It follows that
 \begin{align}
m_{ l, N }(x)=\det &\big\{\sum_{k =1}^{n} c_{1,k}  {\bf p}^{(l)}_{k}e^{\z_kx}, \ldots,
\sum_{k=1}^{n} c_{n,k}  {\bf p}^{(l)}_{k} e^{\z_kx},
\nonumber\\
&\sum_{k =1}^{N} c_{n+1,k}  {\bf p}^{(l)}_{k}e^{\z_kx},\dots, \sum_{k =1}^{N} c_{N-1,k}  {\bf p}^{(l)}_{k}e^{\z_kx}\big\}
\label{eq:mdu}\end{align}
where
\begin{equation}
{\bf p}^{(l)}_{k}=(p_{1,k}, \ldots, p_{l-1,k},p_{l+1,k},\ldots, p_{N,k})^t.
 \label{eq:LP}\end{equation}
 Thus the element with   index $l$ is removed from $ {\bf p} _k =  (p_{1,k}, \ldots,  p_{N,k})^t$ so that in \e{eq:mdu} we take the determinant of
  the $(N-1)\times (N-1)$ matrix. Neglecting columns which are repeated in determinant \e{eq:mdu}, we see that it consists of terms
 \begin{equation}
 \det\big\{   {\bf p}^{(l)}_{1}  e^{\z_1 x} , \ldots, {\bf p}^{(l)}_{n}  e^{\z_n x},
{\bf p}^{(l)}_{k_{n+1}}  e^{\z_{k_{n+1}} x},\ldots, {\bf p}^{(l)}_{k_{N-1}}  e^{\z_{k_{N-1} }x} \big\}
\label{eq:Md}\end{equation}
multiplied by some coefficients which do not depend on $x$. Here the indices $k_{n+1},\ldots, k_{N-1} $ take the values $n+1,\ldots, N $ and $k_p\neq k_q$ if $ p\neq q$.
 Evidently, expression \e{eq:Md} equals
 \begin{equation}
 c_j \exp\big(\tr\!{\bf L}_{0} x -\z_j x\big),\q \tr{\bf L}_{0}= \sum_{k=1}^{N} \z_k,
\label{eq:M5}\end{equation}
for some $j =n+1,\ldots, N$ and a number $c_{j}$ which does not depend on $x$. In view of formulas \e{eq:WW} and \e{eq:U1-}  after division by ${\sf W} (x)$ this expression is $O(| e^{- \z_{j} x}| )$ as $x\to -\infty$.  Hence  \e{eq:as3c} for $j=N$ follows from \e{eq:m2}.
\end{pf}

Let  $\tilde{\bf u}_{j}(x)$ be
   columns   of the matrix $\widetilde{\bf U}(x) $, that is,
   $$
   \widetilde{\bf U}(x)= \{\tilde{\bf u}_{1}  (x),\ldots, \tilde{\bf u}_N  (x)\}.
   $$
        Then  relations \e{eq:as3c} can equivalently be rewritten as
\begin{equation}
\begin{split}
\tilde{\bf u}_{j}  (x)   &= O(e^{\rho_{-} x}),\q x\to -\infty,\q j=n+1, \ldots, N,
\\
\tilde{\bf u}_{j} (x) &= O(e^{-\rho_{+} x}),\q x\to +\infty,\q j= 1, \ldots, n .
\end{split}
\label{eq:as33c}\end{equation}

 Let us define the resolvent kernel $\widetilde{\bf R}(x,y)$ in the same way as in subs.~2.2 with $ {\bf L}(x) $ replaced by $\widetilde{\bf L}(x) $. According to relations  \e{eq:as33c} the first $n$ columns of the  matrix $\widetilde{\bf U}(x) $ exponentially decay as $x\to+\infty$ and its last $N-n$ columns exponentially decay as $x\to-\infty$. Therefore the fundamental matrix $\widetilde{\bf U}(x) $    is  admissible for equation \e{eq:HTd}, and
 it follows from  formula \e{eq:res1} applied to the dual problem  that
 \begin{equation}
\widetilde{\bf R}(x,y)= -\widetilde{\bf U} (x) \widetilde{\sf   P}_{+} \widetilde{\bf U}^{-1} (y) \theta (y-x)
+   \widetilde{\bf U} (x) \widetilde{\sf   P}_{-} \widetilde{\bf U}^{-1} (y) \theta (x-y)
\label{eq:res1du}\end{equation}
where
 \begin{equation}
 \widetilde{\bf U}(x) = {\bf U}^* (x)^{-1}, \q \widetilde{\sf   P}_{\pm}={\sf   P}_{\mp}.
 \label{eq:Udu}\end{equation}
 In particular, we see that
\[
\widetilde{\bf R}(x,y)=- {\bf R}^*(y,x).
\]
This relation implicitly also  follows from the results of subs.~2.3.

 \medskip

{\bf 3.2. }
  Let  $ {\bf G}_{0} (x)=\{g^{(0)}_{j,l}(x)\}$ be the  matrix inverse  to the free matrix \e{eq:Hdc}.  The next result supplements Proposition~\ref{dual} and plays an important role in our proof of the trace formula  \e{eq:tr}.

\begin{proposition}\label{tt}
Let solutions ${\bf u}_{j} (x)$ of equation \e{eq:H2} be distinguished by conditions \e{eq:Solu},  and let $ {\bf U}(x)$ be the corresponding admissible fundamental matrix \e{eq:U}.  Then for all $l=1,\ldots, N$ elements $g_{j,l}  (x) $ of the matrix $ {\bf G}  (x)
={\bf U} ^{-1} (x)$  satisfy the relations
\begin{equation}
\begin{split}
g_{j,l}  (x) -g_{j,l}^{(0)}  (x) &= O(e^{-\rho_{+} x}),\q x\to +\infty,\q j=n+1, \ldots, N,
\\
g_{j,l}  (x) -g_{j,l}^{(0)}  (x)&= O(e^{\rho_{-} x}),\q x\to -\infty,\q j= 1, \ldots, n ,
\end{split}
 \label{eq:as3}\end{equation}
with  positive numbers $\rho_{\pm} $   defined in  \e{eq:min}.
\end{proposition}

\begin{pf}
Let us use the   notation   introduced in the proof of Proposition~\ref{dual}.
We shall again prove the first of relations  \e{eq:as3} for  $j=N$.
According to \e{eq:as2} for $x>\!\!> 0$, we have
 \begin{equation}
m_{ l, N }(x)=\det\big\{\sum_{k =1}^{N} t_{1,k}  {\bf p}^{(l)}_{k}e^{\z_kx}, \ldots,
\sum_{k=1}^{N} t_{n,k}  {\bf p}^{(l)}_{k} e^{\z_kx},
 {\bf p}^{(l)}_{n+1} e^{\z_{n+1} x}, \ldots,  {\bf p}^{(l)}_{N-1} e^{\z_{N-1} x}\big\}
\label{eq:m}\end{equation}
where $ {\bf p}^{(l)}_{k} $ is  vector \e{eq:LP} obtained from $ {\bf p}_{k} $ by removing the component with index $l$.  Neglecting columns which are repeated in $(N-1)\times (N-1)$ matrix \e{eq:m}, we see that
 \begin{align*}
m_{ l, N }(x)=&\det\big\{\sum_{k =1}^{n } t_{1,k}  {\bf p}^{(l)}_{k}  e^{\z_kx}+ t_{1,  N} {\bf p}^{(l)}_{N}    e^{\z_{N}x}, \ldots,
\\
 &\sum_{k=1}^{n }  t_{n,k}  {\bf p}^{(l)}_{k} e^{\z_kx} + t_{n,N}{\bf p}^{(l)}_{N}   e^{\z_{N}x},
{\bf p}^{(l)}_{n+1} e^{\z_{n+1} x}, \ldots, {\bf p}^{(l)}_{N-1} e^{\z_{N-1} x}\big\}  .
\end{align*}
This determinant is the sum of the term
 \begin{equation}
 \det\big\{\sum_{k =1}^{n } t_{1,k}   {\bf p}^{(l)}_{k}  e^{\z_kx}, \ldots, \sum_{k=1}^{n } t_{n,k}   {\bf p}^{(l)}_{k}  e^{\z_kx},
  {\bf p}^{(l)}_{n+1}  e^{\z_{n+1} x}, \ldots,   {\bf p}^{(l)}_{N-1}  e^{\z_{N-1} x}\big\}
\label{eq:m3}\end{equation}
and of the $n$ terms
\[
 \det\big\{t_{1, N}   {\bf p}^{(l)}_{N}  e^{\z_{N} x}, \sum_{k =1}^{n } t_{2,k}  {\bf p}^{(l)}_{k}   e^{\z_kx}, \ldots, \sum_{k=1}^{n } t_{n,k} {\bf p}^{(l)}_{k}  e^{\z_kx},
 {\bf p}^{(l)}_{n+1}  e^{\z_{n+1} x}, \ldots,  {\bf p}^{(l)}_{N-1}  e^{\z_{N-1} x}\big\} ,
\]
\begin{equation}
\cdots
\label{eq:m4}\end{equation}
\[
 \det\big\{ \sum_{k =1}^{n } t_{1,k} {\bf p}^{(l)}_{k} e^{\z_kx}, \ldots,\sum_{k =1}^{n } t_{n-1,k} {\bf p}^{(l)}_{k}  e^{\z_kx} , t_{n,N} {\bf p}^{(l)}_{N} e^{\z_{N} x},
 {\bf p}^{(l)}_{n+1} e^{\z_{n+1} x}, \ldots,  {\bf p}^{(l)}_{N-1} e^{\z_{N-1} x}\big\} .
\]

For the free case where ${\bf V}(x)=0$,  we have the exact equality
\begin{equation}
 m_{ l,N }^{ (0)}(x) =  \det\big\{  {\bf p}^{(l)}_1 e^{\z_1 x}, \ldots,    {\bf p}^{(l)}_n e^{\z_n x},
 {\bf p}^{(l)}_{n+1} e^{\z_{n+1} x}, \ldots,  {\bf p}^{(l)}_{N-1} e^{\z_{N -1} x}\big\}  .
\label{eq:m2Y}\end{equation}
Similarly to \e{eq:U1}, we find that determinant \e{eq:m3} equals
 \begin{equation}
 \det  {\bf T}_{+} \det\big\{  {\bf p}^{(l)}_1 e^{\z_1 x}, \ldots,    {\bf p}^{(l)}_n e^{\z_n x},
 {\bf p}^{(l)}_{n+1} e^{\z_{n+1} x}, \ldots,  {\bf p}^{(l)}_{N-1} e^{\z_{N -1} x}\big\}  .
\label{eq:m6}\end{equation}
By virtue  of  \e{eq:U1} expression \e{eq:m6} divided by  ${\sf W} (x)$ equals
expression \e{eq:m2Y} divided by  ${\sf W}_{0} (x)$.

Thus for the proof of asymptotics  \e{eq:as3} it remains to estimate  determinants \e{eq:m4} by $Ce^{-\rho_{+}x}$. This is similar to the proof  of Proposition~\ref{dual}.
It suffices   to consider  terms corresponding to different values of index $k$ in different sums. Therefore, up to some factors not depending on $x$, determinants \e{eq:m4} consist of the terms
  \begin{equation}
   \det\big\{ {\bf p}^{(l)}_{N} e^{\z_{N  } x}, {\bf p}^{(l)}_{k_{1}}  e^{\z_{k_{1}} x}, \ldots,    {\bf p}^{(l)}_{k_{n-1}}e^{\z_{k_{n-1}} x},
{\bf p}^{(l)}_{n+1} e^{\z_{n+1} x}, \ldots, {\bf p}^{(l)}_{N-1} e^{\z_{N-1  } x}\big\}
\label{eq:m5}\end{equation}
where the indices $k_1,\ldots, k_{n-1} $ take the values $1,\ldots, n $ and $k_p\neq k_q$ if $ p\neq q$. Evidently,   \e{eq:m5} equals expression
\e{eq:M5} for some $j=1,\ldots, n$ and a number $c_j$ which does not depend on $x$. In view of formulas \e{eq:WW} and \e{eq:U1}  after division by ${\sf W} (x)$ this expression is $O(e^{- \rho_{+} x})$ as $x\to +\infty$.
This  concludes the proof of asymptotics  \e{eq:as3}.
 \end{pf}

In view of equality \e{eq:GDd} in terms of  columns $\tilde{\bf u}_{j}  (x) $ of the matrix $\widetilde{\bf U} (x)$  relations \e{eq:as3} can equivalently be rewritten as (cf.  \e{eq:as33c})
\begin{equation}
\begin{split}
\tilde{\bf u}_{j}  (x) -{\bf p}_{j}^* e^{-\bar{\z}_{j}x} &= O(e^{-\rho_{+} x}),\q x\to +\infty,\q j=n+1, \ldots, N,
\\
\tilde{\bf u}_{j} (x) -{\bf p}_{j}^* e^{- \bar{\z}_{j}x}&= O(e^{\rho_{-} x}),\q x\to -\infty,\q j= 1, \ldots, n .
\end{split}
\label{eq:as33}\end{equation}
We emphasize that the functions ${\bf p}_{j}^* e^{-\bar{\z}_{j}x}$    exponentially grow at infinity   while the remainders in formulas  \e{eq:as33} exponentially  decay.

 \medskip

{\bf 3.3. }
Let us, finally, verify that the resolvent kernel defined by equality  \e{eq:res1} satisfies a natural estimate
\begin{equation}
| {\bf R}(x,y)|\leq C e ^{-\rho |x-y|}, \q \rho=\min\{\rho_{+},\rho_{-}\},
\label{eq:RR}\end{equation}
with a  constant\footnote{Here and below $C$   are different positive constants whose precise values are of no importance.} $C$ not depending on $x$ and $y$. It suffices to check that
\begin{equation}
| {\bf U}(x ){\sf P}_{\pm} {\bf U}^{-1}(y) |\leq C e ^{-\rho_{\pm} |x-y|}, \q \pm (y-x)\geq 0.
\label{eq:RR1}\end{equation}
We can suppose that $ {\bf U}(x )$ is defined by formula \e{eq:U} where
 the solutions ${\bf u}_{j}(x)$ of equation  \e{eq:H2}  satisfy condition \e{eq:Solu}.

Note that
\begin{equation}
| {\bf U}(x ){\sf P}_\pm |\leq C e ^{-\rho_\pm |x |}, \q \mp x \geq 0.
\label{eq:RR2}\end{equation}
Passing here to the dual problem and using that $\tilde{\rho}_{\pm}={\rho}_{\mp}$, we see that
\begin{equation}
| {\sf P}_\pm  {\bf U}^{-1}(y ) | = |  {\bf U}^*(y )^{-1}  {\sf P}_\pm  | = | \widetilde{\bf U}(y )\widetilde{\sf P}_\mp     | \leq C e ^{-\rho_\pm |y |}, \q \pm y\geq 0.
\label{eq:RR3}\end{equation}

 Let us prove estimate \e{eq:RR1}, for example, for the upper sign. If $x\leq 0$ and $y\geq 0$, then \e{eq:RR1} directly follows from \e{eq:RR2} and \e{eq:RR3}.

 Suppose next that $x\leq 0$ and $y\leq 0$. According to formula  \e{eq:resZ} it suffices to check that
 \begin{equation}
| u_{k,j}(x) g_{j,l}(y) |  \leq C |e ^{ \z_{j} (x-y)} |
\label{eq:RR4}\end{equation}
for all $k, l= 1,\ldots, N$ and $ j = 1,\ldots, n$. Recall that $| u_{k,j}( x ) | \leq C |e ^{ \z_{j}x } |$ and that
$| g_{j,l}( y ) | \leq C |e ^{ -\z_{j}y } |$ by  Proposition~\ref{tt}. This yields \e{eq:RR4} and hence \e{eq:RR1} for $x\leq 0$,  $y\leq 0$ and the upper sign.

Similarly, we obtain that estimate  \e{eq:RR1} is true  for $x\geq 0$, $y\geq 0$ and the lower sign. Using this estimate for the dual problem we see that
\[
|\widetilde{\bf U}(x )\widetilde{\sf P}_{-} \widetilde{\bf U}^{-1} (y) |\leq C e ^{-\tilde\rho_{-} |x-y|}, \q x\geq 0,\, y\geq 0,\, x\geq y.
\]
 Passing   here to adjoint matrices and taking into account relations  \e{eq:Udu},
  we find that
 \[
| {\bf U}(y ) {\sf P}_{+}  {\bf U}^{-1}(x) |\leq C e ^{- \rho_{+} |x-y|}, \q x\geq 0,\, y\geq 0,\, x\geq y.
\]
Interchanging  $x$ and $y$, we get  \e{eq:RR1} for $x\geq 0$,  $y\geq 0$ and the upper sign.

Thus we have proven the following result.

 \begin{proposition}\label{resest}
Let assumption  \e{eq:zero} hold. Then the   resolvent kernel    \e{eq:res1} satisfies   estimate  \e{eq:RR}.
\end{proposition}

 Estimate  \e{eq:RR} shows that formula \e{eq:resX} obtained for functions $\bf f$ with compact support extends to all ${\bf f}\in L^2({\Bbb R}; {\Bbb C}^N)$.

\section{Trace formula}

 {\bf 4.1. }
 Let us   consider the differential operator \e{eq:H} acting in the space $L^2 ({\Bbb R})$.
  Recall that the operator $H_{0}=i^{-N} \partial^{N}$ is self-adjoint in the space $L^{2}({\Bbb R})$ on domain ${\cal D}(H_{0})$ which is the Sobolev class ${\HH}^N ({\Bbb R})$. The spectrum $\sigma(H_{0})$ of $H_{0}$ coincides with $[0,\infty)$ for $N$ even and with $\Bbb R$ for $N$ odd.
 If the coefficients $v_{j}$, $j= 1,\ldots, N $, of operator \e{eq:V1} belong locally to $L^2$, then the operator
   $H=H_{0}+V$ is well defined by formula  \e{eq:H} at least on the class $C_{0}^\infty ({\Bbb R})$. If all functions $v_{j}$ satisfy  assumption \e{eq:rc}, then according    to Lemma~\ref{RC} the operator $V R_{0}(z)$, $z\not\in \sigma(H_{0})$,  is compact. Therefore the operator $H$ is closed on domain   ${\cal D}(H)={\cal D}(H_{0})$, and by virtue of the Weyl theorem its essential spectrum $\sigma_{ess}(H)=\sigma_{ess}(H_{0})$. The spectrum $\sigma(H)$ of the operator $H$ in
    $ {\Bbb C}\setminus \sigma_{ess}(H_{0})$   consists of eigenvalues (not necessarily real) of finite multiplicity which might accumulate to $\sigma_{ess}(H_{0})$ only.

For the construction of the integral kernel of the resolvent  $R(z)=(H-z)^{-1}$ of the operator $H$, we have to solve the equation
\begin{equation}
 i^{-N} \varphi^{(N)} (x) +  ( V \varphi ) (x)
= z\varphi (x)+ f(x),\q   z \not\in \sigma (H), \q f\in L^2 ({\Bbb R}).
 \label{eq:He}\end{equation}
 Let us rewrite it as a system of differential equations \e{eq:H1}.   We introduce   vectors ${\pmb \varphi}=(\varphi_{1},
  \ldots,   \varphi_{N})^t$ with components $\varphi_{j}=\varphi^{( j-1)}$,  ${\bf f}=i^N(0,0,\ldots, 0, f)^t$ and set
 \begin{equation}
{\bf L}_{0}(z) =
   \begin{pmatrix}
0& 1& \ldots&  0&0
 \\
 0& 0 & \ldots&  0&0
 \\
 \vdots&  \vdots& \ddots&\vdots &  \vdots
  \\
  0&0& \ldots& 0&1
  \\
i^N z &0&  \ldots & 0 &   0
 \end{pmatrix},
\label{eq:J}\end{equation}
\begin{equation}
{\bf V} (x) =  -i^{N }\begin{pmatrix}
0& 0& \ldots&  0&0
 \\
 0& 0 & \ldots&  0&0
 \\
 \vdots&  \vdots& \ddots&\vdots &  \vdots
  \\
  0&0& \ldots& 0&0
  \\
  v_1(x)&     v_{2}(x) &  \ldots &   v_{N-1} (x)&   v_{N }(x)
 \end{pmatrix}.
\label{eq:A}\end{equation}
Then equation \e{eq:He} is equivalent to   vector equation \e{eq:H1}
with the operator ${\bf L}(x,z)$ defined by equality \e{eq:LL}.
We emphasize that it now depends on the spectral parameter $z$.  Matrix \e{eq:J} has eigenvalues $\z_{j}$ such that $\z_{j}^N=i^N z$ with the corresponding eigenvectors
 \begin{equation}
 {\bf  p}_{j} (z) =(1,\z_{j},\ldots, \z_{j}^{N-1})^t.
 \label{eq:pp}\end{equation}
  It is easy to see that
  $n=N/2$ if $N$ is even and $n=(N-1)/2$ for $\Im z>0$ and $n=(N+1)/2$ for $\Im z<0$ if $N$ is odd.

If the coefficients $v_{j} (x)$ have compact supports, then all results of Sections~2 and 3  apply.
 Now we have ${\bf u}_{j}(x,z)=(u_{1,j}(x,z),\ldots, u_{N,j} (x,z))^t$ where the functions $u_{j}(x,z):=u_{1,j}(x,z)$ satisfy homogeneous equation  \e{eq:Heh} and $u_{j}^{(k-1)}(x,z)=u_{k,j}(x,z) $ for $k=2,\ldots, N$. Hence fundamental matrix \e{eq:U} takes form \e{eq:Wrons}. This matrix ${\bf U} (x,z)$ is admissible if the functions $u_{1}(x,z),\ldots, u_{n}(x,z)$ belong to $L^2 ({\Bbb R}_{-})$ and  the functions $u_{n+1}(x,z),\ldots, u_{N}(x,z)$ belong to $L^2 ({\Bbb R}_{+})$.  The Wronskian ${\sf W}(x,z)=\det{\bf U}(x,z)$ satisfies equations \e{eq:Hdet} and \e{eq:WROxm} where
 $\tr {\bf L}(x,z)=-i^{N } v_{N }(x)$. This yields \e{eq:WROx}. Formula \e{eq:resX} means that
  \[
\varphi^{(j-1)} (x) =i^N \int_{-\infty}^\infty R_{j,N} (x,y,z) f(y) dy.
\]
 In particular, we have
  \begin{equation}
  R_{j,N} (x,y,z) =\partial_{x}^{j-1} R_{1,N} (x,y,z),\q j=2, \ldots, N  .
 \label{eq:resY1}\end{equation}
   Condition \e{eq:zero} is equivalent to the assumption that $z$ is not an eigenvalue of the operator $H$. Thus, Proposition~\ref{res} entails the following result.

 \begin{proposition}\label{resH}
Suppose that  $z\not\in\sigma(H)$. Let the matrix ${\bf R} (x,y,z)=\{  R_{j,k} (x,y,z)\}$ be defined by formula \e{eq:res1} where  ${\bf U} (x, z)$ is an admissible fundamental matrix. Then the resolvent $R(z)=(H-z)^{-1}$ of the operator $H$ is the integral operator with kernel
\[
R (x,y,z)=i^N R_{1,N} (x,y,z) .
\]
\end{proposition}

According to formula \e{eq:resZ}  Proposition~\ref{resH} implies that
\begin{equation}
\begin{split}
R (x,y, z)= &- i^N\sum_{j=1}^n  u _{1,j}(x,z) g_{j, N }(y,z),\q  x<y,
\\
R (x,y, z)= &  i^N \sum_{j=n+1}^{N}  u _{1, j}(x,z) g_{j, N}(y,z),\q  x>y,
\end{split}
\label{eq:RESO}\end{equation}

It follows from relations
 \e{eq:Hdc1} and \e{eq:resY1} that, for $N\geq 2$, the function $R (x,y,z)$ as well as its derivatives $R^{(k)}_{x} (x,y,z)$, $k=1,\ldots, N-2$,  are continuous functions of $x$ and $y$ while
\[
R^{(N-1)}_{x} (x,x- 0,z)-R^{(N-1)}_{x} (x,x+0,z)= i^N.
\]

 The case $N=1$ is trivial. Although in this case the kernel $R (x,y, z)$ is not a continuous function, the difference  $R (x,y, z)-R_{0} (x,y, z)$ is continuous and its diagonal values equal zero.

 \medskip

{\bf 4.2. }
Let us find a convenient expression for   integral \e{eq:Di}.
Since $\bf V$ does not depend on $z$,
differentiating equation \e{eq:H3} in $z$, we have
\[
  \dot{\bf U}'  (x,z) = {\bf L} (x,z)  \dot {\bf U}  (x,z)+ \dot {\bf L}_{0} ( z)   {\bf U}  (x,z).
\]
Here and below the dot stands for the derivative in $z$.
Using   formula \e{eq:HT},
we now see that
\begin{align}
  d\big({\bf G}  (x,z)   \dot{\bf U}  (x,z)\big) / dx  = & {\bf G}'  (x,z)   \dot{\bf U}  (x,z)+ {\bf G}  (x,z)   \dot{\bf U}'  (x,z)
  \nonumber\\
  =  & {\bf G}  (x,z) \dot {\bf L}_{0} ( z)  {\bf U}   (x,z)=:  {\bf A}  (x,z).
\label{eq:H3y}\end{align}

Next, we calculate the matrix ${\bf A}  (x,z)$. It follows from  formula  \e{eq:J} that
\[
\dot {\bf L}_{0} ( z)
=  i^N \begin{pmatrix}
0& 0& \ldots&  0
 \\
 \vdots&  \vdots& \ddots&\vdots
 \\
 0& 0& \ldots&  0
  \\
1& 0& \ldots& 0
 \end{pmatrix}
 \]
and hence
\[
{\bf B} (x,z): = \dot {\bf L}_{0} ( z)  {\bf U}  (x,z)
= i^N   \begin{pmatrix}
0& 0& \ldots&  0
 \\
 \vdots&  \vdots& \ddots&\vdots
 \\
 0& 0& \ldots&  0
  \\
u_{ 1,1}& u_{ 1, 2 }& \ldots&  u_{ 1,N}
 \end{pmatrix}
\]
where (cf. \e{eq:Wrons} and \e{eq:U}) $u_{ 1,j}=u_{j}$.
Since ${\bf A}={\bf G}{\bf B}$,
this yields the relation
\begin{equation}
  a_{j,k}  =\sum_{l=1}^{N} g_{j,l}  b_{l ,k} = g_{j,N}  b_{N,k} = i^N  g_{j,N}   u_{ 1, k}
\label{eq:G1}\end{equation}
for the matrix elements $a_{j,k}=a_{j,k}(x,z)$ and $b_{j,k}=b_{j,k}(x,z)$ of the matrices ${\bf A}$ and ${\bf B}$.

Putting together equalities \e{eq:H3y} and \e{eq:G1}, we obtain the relation
\[
i^{N}  u_{1,  k} (x,z) g_{j, N} (x,z) =   \frac{d}{dx} \sum_{l=1}^{N} g_{j,l}(x,z) \dot{u}_{l,k}(x,z)
  \]
  and, in particular,
  \begin{equation}
i^{N}  u_{1, j} (x,z) g_{j, N} (x,z) =   \frac{d}{dx} \sum_{l=1}^{N} g_{j,l}(x,z) \dot{u}_{l,j}(x,z).
  \label{eq:G2}\end{equation}
Comparing formulas \e{eq:RESO} and \e{eq:G2}, we get two representations for diagonal values of the resolvent kernel:
   \[
\begin{split}
R (x,x ,z)  &= -\frac{d}{dx} \sum_{j=1}^n  \sum_{l=1}^N  g_{j, l}(x,z)\dot{u} _{l, j}(x,z),
\\
R (x,x ,  z)  &= \frac{d}{dx} \sum_{j=n+1}^N  \sum_{l=1}^N  g_{j, l}(x,z)\dot{u} _{l, j}(x,z).
\end{split}
\]
Integrating the first of these representations over an  interval $(x_{1}, x)$ and the second over an interval $(x , x_{2})$, we arrive at the following intermediary result.

\begin{proposition}\label{res1}
Under the assumptions of Proposition~\ref{resH} for  all $x_{1}, x_{2}, x \in \Bbb R$, the representation holds:
\begin{equation}
\begin{split}
\int_{x_{1}}^{x_{2} } R (y,y ,  z) dy&=- \sum_{j=1}^{N}  \sum_{l=1}^{N}  g_{j, l}( x ,z)\dot{u} _{l, j}(x ,z)
\\
+ \sum_{j=n+1}^{N}  \sum_{l=1}^{N} & g_{j, l}(x_2,z)\dot{u} _{l, j}(x_{2},z)
+ \sum_{j= 1}^{ n}  \sum_{l=1}^{N}  g_{j, l}(x_{1},z)\dot{u} _{l, j}(x_{1},z).
\end{split}
\label{eq:resint}\end{equation}
\end{proposition}

Let us    consider the first term in the right-hand side of \e{eq:resint}. Since
\[
d \det {\bf U} (x ,z)/dz= \det {\bf U} (x,z) \tr \big( {\bf U} (x, z)^{-1} \dot{\bf U    }(x,z)\big),
\]
we see that
\begin{equation}
\sum_{j=1}^{N}  \sum_{l=1}^{N}  g_{j, l}(x,z)\dot{u} _{l, j}(x,z) =
\tr \big({\bf  G} (x, z) \dot{\bf U  }(x,z)\big)=  {\sf W} (x, z) ^{-1} \dot{\sf W}(x, z)
\label{eq:det1}\end{equation}
where the function $ {\sf W} (x, z) =\det {\bf U} (x ,z)$. Observe that according to
 \e{eq:WROx} this expression does not depend on $x$ which is consistent with formula \e{eq:resint}.

For $N=2$, identity \e{eq:resint} reduces to formula (1.26) of paper  \cite{Finv}. For $N>2$, it is probably new. Identity \e{eq:resint}   allows us to considerably simplify the calculation of $\Tr\big ( R(z) - R_{0}(z)\big)$ compared to the  presentation of  book \cite{Ya} for $N=2$. This  is essential for an arbitrary $N$.

\medskip

{\bf 4.3. }
    Now we suppose that the solutions ${\bf u}_{j}(x,z)$ of equation \e{eq:H2} are distinguished by condition
    \e{eq:Solu} which yields
  condition \e{eq:solu} on the solutions $u_{j}(x,z)$ of equation \e{eq:Heh}; ${\sf W} (x, z)$ is the   Wronskian \e{eq:WRO}.
  Different objects corresponding to the ``free" operator $H_{0}=i^{-N} \partial^N$ will be endowed with index $0$ (upper or lower). For the free case, we   put $u^{(0)}_{j} (x)=e^{\z_{j}x}$, $j=1,\ldots, N$.  Let ${\bf U}_{0}(x,z)$ be the corresponding fundamental matrix \e{eq:Hdc} with the eigenvectors ${\bf p}_{j}(z)$ defined by \e{eq:pp}, and let ${\sf W}_{0}(z)$ be its determinant \e{eq:WRO1}. The normalized Wronskian $\Delta(x , z)$ is defined by formula   \e{eq:det2}.
 We denote by $g_{j, l} (x,z)$ and $g_{j, l}^{(0)}(x,z)$   matrix elements of the matrices  ${\bf G} (x,z)={\bf U}^{-1} (x,z) $ and ${\bf G}_{0}(x,z)={\bf U}_{0}^{-1}(x,z) $, respectively.

     For the proof of  the trace formula,  we combine
Propositions~\ref{tt} and \ref{res1}. Indeed,
let us subtract from equality \e{eq:resint} the same equality for the resolvent $R_{ 0}(z)=(H_{0}-z)^{-1}$ and consider
\begin{equation}
\int_{x_{1}}^{x_{2} } \big( R (y,y ,  z) -R_{0} (y,y ,  z)\big) dy.
\label{eq:resint1}\end{equation}
 Since $\z_{j}^N=i^N z$, for all $l=1,\ldots, N$, we have
\[
\dot{u} _{l, j}(x,z)=\dot{u} _{l, j}^{(0)}(x,z)=d(\z_{j}^{l-1} e^{\z_{j} x})/ dz=i^N N^{-1} \z_{j}^{-N+l-1}(l-1+ x\z_{j} )e^{\z_{j}x}
\]
if either $j=1,\ldots, n$ and $x<\!\!< 0$ or $j=n+1,\ldots, N$ and $x>\!\!> 0$.
According to condition  \e{eq:zeta},  it now directly follows from
 relations   \e{eq:as3}  that,  for all $l=1,\ldots, N$,
 \[
  g_{j, l}(x,z)\dot{u} _{l, j}(x,z)-g^{(0)}_{j, l}(x,z)\dot{u}^{(0)} _{l, j}(x,z) =  O(xe^{-(\rho_{+} + \rho_{-}) |x|})
  \]
  if either $j= 1, \ldots, n$ and $x\to - \infty$ or $ j=n+1,\ldots, N$ and
  $ x\to + \infty$.
 Therefore using equality \e{eq:resint} we see that the contribution to \e{eq:resint1} of
 the terms depending on $x_{1}$ disappear in the limit $x_{1}\to -\infty$ and the contribution  of
 the terms depending on $x_2$ disappear in the limit $x_2\to  +\infty$.
Thus  taking into account equality  \e{eq:det1},
we obtain the following result.

  \begin{theorem}\label{Tr}
Suppose that the coefficients $  v_{1}, \ldots, v_{N }$ of the operator $H$ have compact supports. Then
for   $z\not\in \sigma(H)$, the limit in the left-hand side exists and
\begin{equation}
\lim_{x_{1}\to -\infty,  x_2\to + \infty}\int_{x_{1}}^{x_{2} } \big( R (y,y ,  z) -R_{0} (y,y ,  z)\big) dy  = - \Delta(x, z)^{-1} d\Delta(x, z)/dz
\label{eq:resint2}\end{equation}
where $x\in{\Bbb R}$ is arbitrary.
\end{theorem}

  \medskip

{\bf 4.4. }
  The left-hand side of \e{eq:resint2} can of course be identified with   the trace of the difference   $ R(z)-R_{0}(z) $. To that end, we first   verify inclusion  \e{eq:TSS}. We proceed from
   the following
 well-known result (see, e.g., survey \cite{BS2} by M.~Sh.~Birman and M.~Z.~Solomyak).

   \begin{proposition}\label{TrCl}
Suppose that
\begin{equation}
{\cal B}^2:=\int_{-\infty}^{\infty } (1+\xi^2)^{\alpha} |b(\xi)|^2 d\xi <\infty , \q
{\cal V}^2:=  \int_{-\infty}^{\infty } (1+x ^2)^{\alpha} |v(x)|^2 dx <\infty
 \label{eq:TrCl}\end{equation}
for some $\alpha> 1/2$.  Then the integral operator $T: L^2 ({\Bbb R}; dx)\to L^2 ({\Bbb R}; d\xi)$ with kernel \e{eq:TrCl2}
belongs to the trace class ${\goth S}_{1}$ and its trace norm satisfies the estimate
$\| T\| _{{\goth S}_{1}}  \leq C {\cal B}\, {\cal V} $
 where the constant $C$ depends on $\alpha>1/2$ only.
\end{proposition}

Now it is easy to prove the following result.

   \begin{lemma}\label{TrC}
Under assumption  \e{eq:V}  inclusion  \e{eq:TSS} holds.
\end{lemma}

  \begin{pf}
  Let us proceed from the resolvent identity
 \begin{equation}
R(z)-R_{0}(z)= -R_{0}(z) V R (z)=-\sum_{j=1}^{N } R_{0}(z) v_{j}\partial^{j-1} R (z),\q z\not\in\sigma(H),
 \label{eq:ResId}\end{equation}
 where $v_{j}$ is the operator of multiplication by the function $v_{j} (x)$. The operators $\partial^j R (z)$ are bounded because the operator $H_{0}  R (z)$ is bounded.
 Proposition~\ref{TrCl} implies that $R_{0}(z) v_{j}\in {\goth S}_{1}$  if $N>1$. Indeed, let $\Phi: L^2 ({\Bbb R}; dx)\to L^2 ({\Bbb R}; d\xi)$ be the Fourier transform. Then the operator $\Phi R_{0}(z) v_{j}$ has integral kernel  \e{eq:TrCl2} where $b(\xi) =(2\pi)^{-1/2} (\xi^N-z)^{-1}$ and $v(x)=v_{j}(x)$. By virtue of
 \e{eq:V} condition \e{eq:TrCl} is satisfied in this case. If $N=1$, then we can use that the operator $|v_{1}|^{1/2} R_{0}(z)$ belongs to the Hilbert-Schmidt class. Thus all terms in the right-hand side of \e{eq:ResId} belong to the trace class.
   \end{pf}

 Since the kernel $ R (x,y ,  z)$ (and $ R_{0} (x,y ,  z)$) is a continuous function of $x,y$ and the limit in the left-hand side of \e{eq:resint2} exists, we see (see, e.g., \cite{Ya}, Proposition~3.1.6) that
\begin{equation}
\Tr \big( R(z)-R_{0}(z)\big)=
 \int_{-\infty}^{\infty } \big( R (y,y ,  z) -R_{0} (y,y ,  z)\big) dy   .
\label{eq:resint2w}\end{equation}
Putting together formulas \e{eq:resint2} and \e{eq:resint2w}, we get the trace formula  \e{eq:tr} for the   coefficients $  v_{1}, \ldots, v_{N}$ with compact supports.

\section{Integral equations}

  Here we consider differential equation  \e{eq:Heh} with arbitrary short-range coefficients.   Actually, we follow the scheme of Section~2 and first consider a more general  equation  \e{eq:H2}.

  \medskip

 {\bf 5.1.}
 As usual, we suppose that the eigenvalues $\z_{j}$, $j=1,\ldots, N$,  of an $N\times N$ matrix ${\bf L}_{0}$ are distinct and do not lie on  the imaginary axis.  We set
  \begin{equation}
  {\bf P}_{j}= \langle \cdot,  {\bf p}_{j}^*\rangle   {\bf p}_{j}, \q j=1,\ldots, N,
  \label{eq:Pp}\end{equation}
  where $ {\bf p}_{j}$ are eigenvectors of ${\bf L}_{0}$  and the vectors $ {\bf p}_{j}^*$
form  the dual basis. We have ${\bf P}_{j}^2= {\bf P}_{j}$, ${\bf P}_{j} {\bf P}_{k}=0$ if $j\neq k$, and
\begin{equation}
{\bf L}_0{\bf P}_{j}= \z_{j} {\bf P}_{j}, \q\sum_{j=1}^N {\bf P}_{j}={\bf I}.
\label{eq:di}\end{equation}

   Let a matrix ${\bf L}(x)$ be given by formula \e{eq:LL} where we now assume that
    \begin{equation}
{\bf V}\in L^1 ({\Bbb R_{\pm}}).
\label{eq:VV}\end{equation}
We shall show that, for all $j=1,\ldots, N$, equation  \e{eq:H2} has solutions ${\bf u}_{j}^{(\pm)}(x)$ such that
\begin{equation}
{\bf u}_{j}^{(\pm)}(x)=  e^{\z_{j} x}( {\bf p}_{j}  +o(1)),\q x\to \pm\infty.
\label{eq:VV1}\end{equation}
Thus we construct solutions of \e{eq:H2} both (exponentially) decaying and growing at infinity.
We emphasize that our construction of the functions ${\bf u}_{j}^{(+)}(x)$ (of ${\bf u}_{j}^{(-)}(x)$) requires   condition \e{eq:VV} for $x\in {\Bbb R}_{+}$ (for $x\in {\Bbb R}_-$) only.  Functions ${\bf u}_{j}^{(\pm)}(x)$  will be defined as solutions of integral equations which we borrow from the book \cite{CL} (see Problem~29 of Chapter~3). For definiteness, we consider the case $x\to -\infty$ and put ${\bf u}_{j} ={\bf u}_{j}^{(-)}$.

 Set
  \begin{equation}
{\bf K}_{j}(x )=  \sum_{m: \kappa_{m}> \kappa_{j}}{\bf P}_m  e^{{\z}_m  x  }\theta(-x)
- \sum_{m: \kappa_{m} \leq\kappa_{j}}  {\bf P}_m  e^{{\z}_m  x  }\theta(x),
\label{eq:IE1a}\end{equation}
where $\kappa_m=\Re \z_m$. It follows from relations \e{eq:di} that
  \begin{equation}
{\bf K}_{j}'(x )= {\bf L}_{0} {\bf K}_{j}(x )-\delta (x){\bf I} ,
\label{eq:IE1d}\end{equation}
where $\delta (x)$ is the Dirac function.
We also use the estimate
 \begin{equation}
|{\bf K}_{j}(x )|\leq C_{j} \big(\sum_{m: \kappa_{m}> \kappa_{j}}   e^{\kappa_m  x  }\theta(-x) + e^{\kappa_{j} x} \theta(x) \big),
\label{eq:Gest1}\end{equation}
which is a direct consequence of definition \e{eq:IE1a}. In particular, we see that
  \begin{equation}
|{\bf K}_{j}(x )|\leq C_{j} e^{\kappa_{j} x}, \q \forall x\in {\Bbb R}.
\label{eq:Gest}\end{equation}

 Let $\chi_{X}$ be the characteristic function of an interval $X$.
Consider the integral equation
 \begin{equation}
 {\bf u}_{j}(x  )=  e^{{\z}_{j} x} {\bf p}_{j}
- \int_{-\infty}^{\infty} {\bf K}_{j} (x-y ) {\bf V} (y  ) \chi_{(-\infty, a)}(y)  {\bf u}_{j}(y  ) dy,\q x<a,
\label{eq:IE1}\end{equation}
for a function $ {\bf u}_{j}(x  )= {\bf u}_{j}(x ;a )$ depending on    the parameter $a$ which will be chosen later.   If $\kappa_{j}=\max_{m}\kappa_m$, then the first sum in \e{eq:IE1a} is absent. In this case we can omit $\chi_{(-\infty, a)}(y)$ in \e{eq:IE1} so that  \e{eq:IE1} becomes a Volterra integral equation.   However \e{eq:IE1} is only a Fredholm equation for other values of $j $.  Suppose that a function ${\bf u}_{j}(x  )$ satisfies the estimate $ {\bf u}_{j}(x  )=O(e^{\kappa_{j} x} )$ as $x\to-\infty$ and   equation \e{eq:IE1}.
Differentiating \e{eq:IE1} and using \e{eq:IE1d} we see that
\begin{equation}
 {\bf u}'_{j}(x  )=  {\z}_{j} e^{{\z}_{j} x} {\bf p}_{j}
+ {\bf V} (x  ) \chi_{(-\infty, a)}(x)  {\bf u}_{j}( x  )
- {\bf L}_{0}\int_{-\infty}^{\infty} {\bf K}_{j} (x-y ) {\bf V} (y  ) \chi_{(-\infty, a)}(y)   {\bf u}_{j}(y  ) dy .
\label{eq:IE1der}\end{equation}
Putting together equations \e{eq:IE1} and \e{eq:IE1der} we find that a solution  ${\bf u}_{j}(x  )$ of integral equation  \e{eq:IE1} satisfies also the differential equation
\[
 {\bf u}'_{j}(x  )= {\bf L}_{0} {\bf u}_{j}(x  )+ {\bf V} (x  ) \chi_{(-\infty, a)}(x)  {\bf u}_{j}(x  ),
\]
which reduces to  equation \e{eq:H2} for $x<a$.

   Let us set
    \begin{equation}
 {\bf u}_{j}(x ;a)= e^{{\z}_{j} x}{\bf w}_{j}(x;a),\q x<a,
\label{eq:Eq}\end{equation}
and rewrite equation \e{eq:IE1} as
 \begin{equation}
{\bf w}_{j}(x ;a )= {\bf p}_{j}
- \int_{-\infty}^{a } {\bf K}_{j} (x-y )e^{-{\z}_{j} (x-y)}  {\bf V} (y  )   {\bf w}_{j}(y ;a ) dy .
\label{eq:IE1w}\end{equation}
By virtue of    assumption \e{eq:VV} and estimate \e{eq:Gest}  we   can choose   the parameter $a $ such that
  \begin{equation}
\int_{-\infty}^{a }\big| {\bf K}_{j} (x-y )e^{-{\z}_{j} (x-y)}  {\bf V} (y  )\big|  dy \leq
C \int_{-\infty}^{a }\big|   {\bf V} (y  )\big|  dy <1,\q \forall x\in{\Bbb R}.
\label{eq:Eq1s}\end{equation}
Then   equation \e{eq:IE1w} can   be solved in the space $L^\infty ((-\infty,a); {\Bbb C}^N)$ by the method of successive approximations.

This result   can also be reformulated in the following way. Let  ${\Bbb Q}_{j} (a)$ be the   integral operator with   kernel
\begin{equation}
{\bf Q}_{j} (x,y )= {\bf K}_{j} (x-y )e^{-{\z}_{j} (x-y)}  {\bf V} (y  )
\label{eq:QQ}\end{equation}
acting  in the space $L^\infty ( (-\infty,a) ; {\Bbb C}^N )$.
 Then
\begin{equation}
{\bf w}_{j}( a )= \big(I- {\Bbb Q}_{j} (a))^{-1}{\bf p}_{j}
\label{eq:IE1Q}\end{equation}
where the inverse operator exists because  $\| {\Bbb Q}_{j} (a)\| <1$.

Clearly, the function ${\bf u}_{j} (x;a )$ defined by formula   \e{eq:Eq} satisfies differential equation \e{eq:H2}  for $x<a$. Since a solution of a differential equation of first order is determined uniquely by its value at one point, it suffices to require equality \e{eq:Eq} only for one $x<a$. Then the corresponding solution can be extended to all $x\in \Bbb R$. Now we are in a position  to give the precise definition.

\begin{definition}\label{DUa}
 Let     ${\bf w}_{j}^{(-)}(\cdot;a_{-})\in L^\infty ((-\infty,a_{-}); {\Bbb C}^N)$, $j=1,\dots,N$,  be the  solution of equation \e{eq:IE1w} where
 $a=a_{-}$ is a sufficiently large negative number. We define ${\bf u}_{j}^{(-)}(x;a_{-})$ as   the solution of differential equation \e{eq:H2} which satisfies   condition \e{eq:Eq} for some (and then for all) $x <a  $. The solutions ${\bf u}_{j}^{(+)}(x;a_{+})$, $j=1,\dots,N$,   are defined quite similarly if $(-\infty,a_{-})$ is replaced by $(a_{+}, \infty)$  where
 $a_{+}$ is a sufficiently large positive number.
 \end{definition}

  It remains to verify asymptotics \e{eq:VV1} for the function $ {\bf u}_{j}(x )$. According to  \e{eq:Eq} and  \e{eq:IE1w}   it suffices to check that the integral in the right-hand side of  \e{eq:IE1w} tends to zero as $x\to-\infty$. Using estimate \e{eq:Gest1} and the inclusion  ${\bf w}_{j}\in L^\infty ((-\infty,a_{-}); {\Bbb C}^N )$, we see that this integral is bounded by
  \begin{equation}
C \int_{-\infty}^{x}    \big| {\bf V} (y  )\big| dy + C \sum_{m: \kappa_{m}> \kappa_{j}}\int_{x}^{a}     e^{(\kappa_{m}- \kappa_{j})  (x-y)  } \big| {\bf V} (y  )\big| dy.
\label{eq:Eq1ss}\end{equation}
The first integral here tends to zero as $x\to-\infty$ by virtue of condition \e{eq:VV}. Each of the integrals over $(x, a )$ can be estimated by
\[
e^{(\kappa_{m}- \kappa_{j}) x/2  } \int_{x/2}^{a }       \big| {\bf V} (y  )\big| dy
+\int_{x}^{x/2}       \big| {\bf V} (y  )\big| dy.
\]
Since $\kappa_{m}> \kappa_{j}$, this expression tends to zero as $x\to-\infty$ by virtue again of condition \e{eq:VV}.

  Thus we arrive at the following result.

 \begin{proposition}\label{SRD}
 Let   assumption \e{eq:VV} hold, and let $a_+$ $(a_-)$ be a sufficiently large positive $($negative$)$ number.  Then, for all   $j=1,\ldots, N$, the functions
 $ {\bf u}_{j}^{(\pm)}(x; a_{\pm} )$  $($see Definition~\ref{DUa}$)$ satisfy  equation \e{eq:H2} and have asymptotics \e{eq:VV1}.
\end{proposition}

 Solutions   of    equation \e{eq:H2} are  of course not  determined uniquely by asymptotics \e{eq:VV1}.
    In particular,   the solutions  ${\bf u}_{j}^{(\pm)}(x; a_{\pm} )$ generically depend on the choice of the parameter $a_{\pm} $.

    Let   ${\bf u}_{j} (x; a, r ) = {\bf u}_{j}^{(\pm)}(x; a_{\pm},r )$ be the function constructed above
    for the cut-off coefficient ${\bf V}_{r}  (x )=\chi_{(-r,r)}(x)  {\bf V}  (x )$; thus function    \e{eq:IE1Q}   is now replaced by
    \begin{equation}
{\bf w}_{j}( a,r )= \big(I- {\Bbb Q}_{j} (a) \chi_{(-r,r)} )^{-1}{\bf p}_{j}.
\label{eq:IE1Qr}\end{equation}
        Since
        $$
     \lim_{r\to \infty}   \|{\Bbb Q}_{j} (a) \big(1- \chi_{(-r,r)} \big) \|_{L^\infty (-\infty, a)} = 0,
        $$
        we see that
        $ {\bf u}_{j}  (x;a, r ) \to {\bf u}_{j}   (x ,a )$ as $r\to\infty$
         for all fixed $x<a $.
  This relation extends to all $x\in{\Bbb R}$ because    solutions of differential equations depend continuously on initial data.
Therefore   Proposition~\ref{SRD} can be supplemented by the following result.

  \begin{lemma}\label{SRD2}
 Under the  assumptions of Proposition~\ref{SRD}, let ${\bf u}_{j}^{(\pm)}  (x; a_{\pm}  )$ and ${\bf u}_{j}^{(\pm)}  (x; a_{\pm},  r )$ be the solutions of equations \e{eq:H2}  with  $ {\bf V}  (x )$ and $  {\bf V}_{r}  (x )$, respectively, specified in Definition~\ref{DUa}. Then
 for all $j=1,\ldots, N$ the  relation
  \begin{equation}
\lim_{r\to\infty}{\bf u}_{j}^{(\pm)}  (x; a_{\pm} , r ) = {\bf u}_{j}^{(\pm)}  (x ; a_{\pm}  )
 \label{eq:lim}\end{equation}
 holds  uniformly in $x$ on compact intervals of $\Bbb R$.
 \end{lemma}

  \medskip

  {\bf 5.2.}
    If   a function ${\bf u}_{j}^{(\pm)}(x )$ satisfies   equation \e{eq:H2} and has asymptotics \e{eq:VV1}, then   adding to ${\bf u}_{j}^{(\pm)}(x )$ a solution with a more rapid decay (or less rapid growth) as $x\to\pm \infty$ we obtain again a  solution   of
    equation \e{eq:H2} with the same asymptotics \e{eq:VV1}.
  It is natural to expect that  this procedure exhausts the  arbitrariness in the definition of
  ${\bf u}_{j}^{(\pm)}(x )$.  The precise result will be formulated in Lemma~\ref{LIx}.

    The following assertion is almost obvious.

\begin{lemma}\label{LI}
Suppose that solutions ${\bf u}_{1}^{(\pm)},\ldots, {\bf u}_N^{(\pm)}$ of the differential equation \e{eq:H2} have asymptotics \e{eq:VV1}  as   $x\to \pm\infty$.
   Then for each of the signs $``\pm"$ the functions ${\bf u}_{1}^{(\pm)},\ldots, {\bf u}_N^{(\pm)}$ are linearly independent.
   \end{lemma}

 \begin{pf}
 It follows from \e{eq:VV1} that
 \[
 \det\{{\bf u}_{1}^{(\pm)} (x),\ldots, {\bf u}_N^{(\pm) } (x) \}
 =\det\{{\bf p}_{1} ,\ldots, {\bf p}_N \}\exp\big(\sum_{j=1}^{N}\z_{j} x \big)
 \big(1+ o(1)\big)
 \]
 as   $x\to \pm\infty$. Since $\det\{{\bf p}_{1} ,\ldots, {\bf p}_N \}\neq 0$,
  this  expression   is not zero for sufficiently large $\pm x$.
   \end{pf}

   \begin{lemma}\label{LIx}
Suppose that solutions ${\bf u}_{1}^{(\pm)},\ldots, {\bf u}_N^{(\pm)}$ of  the differential equation \e{eq:H2} have asymptotics \e{eq:VV1}  as   $x\to \pm\infty$. Let
$\tilde{\bf u}_j^{(\pm)} $ be an arbitrary  solution of   \e{eq:H2} with asymptotics \e{eq:VV1}  as   $x\to \pm\infty$.
   Then necessarily
    \begin{equation}
\tilde{\bf u}_j^{(\pm)} (x) ={\bf u}_j^{(\pm)} (x)+
\sum_{\pm (\kappa_{l}-\kappa_{j})< 0} c_{j,l}^{(\pm)} {\bf u}_l ^{(\pm)} (x)
 \label{eq:Wro1}\end{equation}
 for some numbers $c_{j,l}^{(\pm)}$.
   \end{lemma}

 \begin{pf}
As before, we set
    ${\bf u}_{j} (x)={\bf u}_{j}^{(-)}(x)$.      According to Lemma~\ref{LI} we have
     \begin{equation}
\tilde{\bf u}_j  (x) =
\sum_{l=1}^N c_{j,l}  {\bf u}_l   (x)
 \label{eq:Wro2}\end{equation}
 with some numbers $c_{j,l}$.
 Therefore it follows from asymptotic relations \e{eq:VV1} that
  \[
  e^{\z_{j}x}({\bf p}_j  +o(1)) =
\sum_{\kappa_{l}\leq \kappa_{j}}  c_{j,l}   e^{\z_l x}({\bf p}_l +o(1))+ o(e^{\kappa_{j}  x})
\]
  as $x\to-\infty$.
  Since the vectors ${\bf p}_1,\dots, {\bf p}_N$ are linearly independent, this relation implies that $c_{j,l}=0$ if $l\neq j$ and $c_{j,j}=1$.
  Thus equality \e{eq:Wro2} leads to \e{eq:Wro1}.
  \end{pf}

    \medskip

  {\bf 5.3.}
  Let us return to differential equation  \e{eq:Heh} with   coefficients satisfying the assumption
   \begin{equation}
v_{j}\in L^1 ({\Bbb R}_{\pm}), \q j=  1, \ldots, N ,
\label{eq:Vv}\end{equation}
 only.
 Define as usual the matrices ${\bf L}_{0}(z)$ and ${\bf V} (x)$   by formulas \e{eq:J} and \e{eq:A}. Now $\z_{j}^N= i^N z$ and the vectors ${\bf p}_{j} (z)$, $j=1,\ldots, N$, are given by formula \e{eq:pp}. It is easy to control the dependence on $z$ of matrices \e{eq:Pp}.

 \begin{lemma}\label{Pp11}
 Elements   $p^{(j)}_{k,l}(z)$ of the matrices ${\bf P}_{j}(z)$, $j=1,\ldots, N$,  obey the relation
 \begin{equation}
 p^{(j)}_{k,l}(z)
  = O (| \z|^{k-l}), \q |\z|^N= |z| \to\infty.
 \label{eq:Pp11}\end{equation}
\end{lemma}

 \begin{pf}
 Obviously, the basis dual to ${\bf p}_{j} (z)$ consists of the vectors
 $$
 {\bf p}_{j}^* (z)= (c_{j,1}, c_{j,2}\z_{j}^{-1},\ldots, c_{j,N}\z_{j}^{-N+1})
 $$
  where the numbers $ c_{j,l}$ do not depend on $|z|$. It now follows from equality  \e{eq:pp} that
   $
   p^{(j)}_{k,l}(z)= \bar{c}_{j,l}\z_{j}^{k-1}\bar{\z}_{j}^{-l+1},
    $
    which yields \e{eq:Pp11}.
   \end{pf}

  The next step is to control the dependence on $z$ of matrices \e{eq:QQ}.

 \begin{lemma}\label{Pp12}
 Elements   $q^{(j)}_{k,l}(x,y,z)$ of the matrices ${\bf Q}_{j}(x,y,z)$, $j=1,\ldots, N$,
 admit the  estimate
   \begin{equation}
| q_{k,l}^{(j)}(x,y, z)| \leq
C | \z|^{k-N} | v_l (y)|
  \label{eq:ele}\end{equation}
    with a constant $C$ not depending on $x$, $y$ and $z$.
\end{lemma}

 \begin{pf}
 It follows from \e{eq:IE1a} and \e{eq:Pp11} that  elements $s^{(j)}_{k,l}(x,   z)$ of  the matrix $ {\bf K}_{j} (x )e^{-\z_{j} x} $ satisfy the estimate (cf. \e{eq:Gest})
   \[
| s_{k,l}^{(j)}(x,z)| \leq
C | \z|^{k-l}.
\]
  This directly implies \e{eq:ele} because, by definition \e{eq:A},
  $$
  q_{k,l}^{(j)}(x,y, z)= -i^N s_{k,N}^{(j)}(x-y,z)   v_l (y).
  $$
  \end{pf}

 Let ${\bf u}_{j}^{(\pm)}(x,z ;a_{\pm})$ (recall that ${\bf u}_{j}^{(\pm)}=(u_{1,j}^{(\pm)},\ldots, u_{N,j}^{(\pm)})^t$) be the solution of equation \e{eq:H2} specified in Definition~\ref{DUa}. Then the function
  \begin{equation}
  u_{j}^{(\pm)}(x,z;a_{\pm} ):= u_{1, j}^{(\pm)}(x,z;a_{\pm})
  \label{eq:Qq}\end{equation}
   satisfies also equation \e{eq:Heh} and according to \e{eq:pp} asymptotics \e{eq:VV1} imply asymptotics \e{eq:as10}. Therefore  Proposition~\ref{SRD} and Lemma~\ref{SRD2} entail the following result.

 \begin{proposition}\label{SRD1}
 Let   assumption \e{eq:Vv} hold, let $|z|\geq c>0$ and let $a_{+}=a_{+}(c)$ $(a_{-}=a_{-}(c))$ be a sufficiently large positive $($negative$)$ number.  Then
 for every $j=1,\ldots, N$  the function $u_{j}^{(\pm)}(x, z;a_{\pm} )$ determined by Definition~\ref{DUa} and equality \e{eq:Qq} satisfies equation \e{eq:Heh} and has  asymptotics \e{eq:as10} as $x\to\pm \infty$.
  Moreover,  the corresponding solutions $u_{j}^{(\pm)}(x,z;a_{\pm}, r  )$ of equation \e{eq:Heh} with
  cut-off coefficients $\chi_{(-r,r)}(x)  v_{k}  (x )$, $k= 1,\ldots, N $, satisfy the relation
   \begin{equation}
\lim_{r\to\infty}u_{j}^{(\pm)}  (x,z; a_{\pm}, r ) = u_{j}^{(\pm)}  (x,z; a_{\pm}  )
\label{eq:cut2}\end{equation}
 uniformly in $x$ on compact intervals of $\Bbb R$.  This relation remains true for $N-1$ derivatives of the functions $u_{j}^{(\pm)}$.
\end{proposition}

By definition  \e{eq:IE1a}, the   kernels ${\bf K}_{j}(x,z)$  depend analytically on $z$ except on the rays where $\Re\z_l = \Re\z_{j}$   for  some root $\z_l \neq \z_{j}$   of the equation $\z ^N= i^N z$. In addition to the condition $z\not\in\sigma(H_{0})$, this excludes also the half-line $z<0$ for even $N$ and the line $\Re z=0$ for odd $N$.  Hence the same
is true for the functions $u_{j}^{(\pm)}(x, z;a_{\pm} )$ if $|z|>c>0$. Thus, for fixed $x$ and $a_{\pm}$, the functions $u_{j}^{(\pm)}(x, z;a_{\pm} )$  are analytic functions of $z$ if $|z| > c>0$, $\Im z\neq 0$ for $N$ even and  if $\Im z\neq 0$, $\Re z\neq 0$  for $N$ odd. On the rays where $\Re\z_l =\Re\z_{j}$, the limits of  $u_{j}^{(\pm)}(x, z;a_{\pm} )$ from both sides exist but differ, in general, from each other by a term which decays faster (or grows less rapidly) than $e^{\z_{j} x}$ as $x\to+ \infty$ or
$x\to - \infty$.

\medskip

  {\bf 5.4.}
  Integral equations \e{eq:IE1} turn out   also to be useful (even for functions $v_{j}(x)$ with compact supports) for a study of asymptotics of the solutions $u_{j}^{(\pm)}(x,z;a_{\pm} )$  of differential equation  \e{eq:Heh}  as $| z | \to \infty$.
  We choose the sign $``-"$, fix the parameter $a=a_{-}$ and index $j$ and drop them out of notation.

          Consider   system \e{eq:IE1w} of $N$ equations for components $w_{k}  (x,z)$, $k=1,\ldots, N$, of a vector-valued function ${\bf w} (x,z)$.
      Set $w_{k}  (x,z)=\z^{k-1} \tilde{w}_{k}  (x,z)$ and take into account equality  \e{eq:pp}. Then we obtain for $\tilde{w}_{k}  (x,z)$ a system
    \begin{equation}
\tilde{w}_{k}  (x,z)= 1- \sum_{l=1}^{N } \z^{l-k} \int_{-\infty}^a q_{k,l} (x, y,z)  \tilde{w}_l (y,z) dy , \q x< a,
   \label{eq:eleme}\end{equation}
   where the elements $q_{k,l} $ of the matrix ${\bf Q}$ satisfy inequality \e{eq:ele}.
   In particular, the operator in the right-hand side of  \e{eq:eleme} is uniformly bounded
as   $|z|\to \infty$.

   Assume additionally that  $v_{N }(x)=0$. Then according to \e{eq:ele}   the norm of the operator in the right-hand side of  \e{eq:eleme} is $ O (  |\z|^{-1})$ as $|z|\to \infty$. Therefore for sufficiently large $|\z|$ system \e{eq:eleme} can be solved in the space $L^\infty (( -\infty,a) ; {\Bbb C}^N )$ by the method of successive approximations and $\tilde{w}_{k}  (x,z)=1+ O (  |\z|^{-1})$, $k=1,\ldots, N$.

 As we have already seen in the proof of Proposition~\ref{SRD}, the solution of system  \e{eq:eleme} necessarily has the asymptotics $\tilde{w}_{k}  (x,z)=1+o(1)$, $k=1,\ldots, N$,  as $x\to -\infty$. Define as usual $u(x,z)$ as a solution of equation \e{eq:Heh} such that $u(x,z)=e^{\z x} \tilde{w}_1 (x,z)$ for $x<a$. Thus we obtain   the following result.

 \begin{proposition}\label{HE}
  Let   assumption \e{eq:Vv} hold, and let $v_{N}=0$. Fix arbitrary $a_{\pm}$.
    Then
 for all $j=1,\ldots, N$ and all  sufficiently large $|z|$   equation \e{eq:Heh} has   solutions $u_{j}^{(\pm)}(x , z;  a_{\pm})$ with    asymptotics \e{eq:as10} as $x\to\pm\infty$ and such that
 \begin{equation}
u_{j}^{(\pm)}(x,z;  a_{\pm} )=   e^{\z_{j} x}\big( 1+ O (  |z|^{-1/N})\big), \q | z | \to \infty,
\label{eq:Lx}\end{equation}
for  all $\pm (x -a_{\pm})>0$.
\end{proposition}

 \begin{remark}\label{HEr}
If   assumption \e{eq:Vv} is true for both signs, then we can set $a=+\infty$ in equation \e{eq:eleme}. For sufficiently large $|z|$, such an equation can again be solved by the method of successive approximations.
\end{remark}

\section{The Wronskian and the perturbation determinant}

 {\bf  6.1.}
   Let us  define the Wronskian ${\sf W}(x)$ for differential equation   \e{eq:H2}
   where the matrix-valued function ${\bf L} (x)$ is given by formula \e{eq:LL} and
${\bf V} (x)$     satisfies assumption \e{eq:VV} (for both signs) only.   To justify the definition below, we start with the following observation.

        \begin{lemma}\label{LI1}
Suppose that both sets of solutions ${\bf u}_{j} $ and  $\tilde{\bf u}_{j} $
  of the differential equation \e{eq:H2} have asymptotics \e{eq:VV1} as $x\to -\infty$ for $j=1,\ldots, n$  and  as $x\to +\infty$   for $j=n+1,\ldots, N$.
   Then for all $x$
    \begin{equation}
 \det
\{  {\bf u}_{1} (x),\ldots,  {\bf u}_{N}(x ) \}
=\det \{  \tilde{\bf u}_{1} (x ), \ldots, \tilde{\bf u}_{N}(x ) \} .
\label{eq:LI2}\end{equation}
   \end{lemma}

    \begin{pf}
Let us proceed from    Lemma~\ref{LIx}. To simplify notation, we suppose that
  \begin{equation}
  \kappa_{1}\geq\ldots\geq \kappa_{n}>\kappa_{n+1}\geq\ldots\geq \kappa_N.
  \label{eq:or}\end{equation}
First, we check that for all $l=1, \ldots, n$,
 \begin{equation}
\det \{  \tilde{\bf u}_{1} (x ), \ldots, \tilde{\bf u}_{N}(x ) \}
=\det \{   {\bf u}_{1} (x ), \ldots,  {\bf u}_{l} (x ),  \tilde{\bf u}_{l+1} (x ),  \ldots, \tilde{\bf u}_{N}(x ) \} .
\label{eq:LIa}\end{equation}
For $l=1$ this equality is obvious because necessarily $\tilde{\bf u}_{1} (x )= {\bf u}_{1} (x )$. Suppose that \e{eq:LIa} is true for some $l$. Then using \e{eq:Wro1} we see that the left-hand side of \e{eq:LIa} equals
 \[
\det \{   {\bf u}_{1} (x ), \ldots,  {\bf u}_{l} (x ),  {\bf u}_{l+1} (x ) +\sum_{\kappa_m>\kappa_{l+1}}c_{l+1,m}{\bf u}_m (x ), \tilde{\bf u}_{l+2} (x ), \ldots,  \tilde{\bf u}_{N}(x ) \} .
\]
Since according to \e{eq:or} the contribution of the sum over $\kappa_{m}>\kappa_{l+1}$ equals zero, this yields
relation  \e{eq:LIa} for $l+1$ and hence for all $l=1, \ldots, n$.

Quite similarly, we can verify that for all $l=   N, \ldots, n+1$
 \begin{align}
\det \{ & {\bf u}_{1} (x ), \ldots,  {\bf u}_{n} (x ), \tilde{\bf u}_{n+1} (x ), \ldots, \tilde{\bf u}_{N}(x ) \}
\nonumber\\
&=  \det \{  {\bf u}_{1} (x ), \ldots,  {\bf u}_{n} (x ), \tilde{\bf u}_{n+1} (x ), \ldots, \tilde{\bf u}_{l}(x ),  {\bf u}_{l+1} (x ), \ldots,  {\bf u}_{N} (x ) \}.
\label{eq:LIb} \end{align}

Putting together equalities  \e{eq:LIa}  and \e{eq:LIb} for $l=n$, we get \e{eq:LI2}.
   \end{pf}

   Now we are in a position to define the Wronskian $ {\sf W} (x )$.

      \begin{definition}\label{WRoo}
 Let ${\bf u}_j   (x)$, $j=1,\ldots, N$,   be {\it arbitrary} solutions of equation \e{eq:H2} with asymptotics \e{eq:VV1} as $x\to -\infty$ if $j=1,\ldots, n$ and  as $x\to \infty$ if $j=n+1,\ldots, N$. We set
  \begin{equation}
 {\sf W} (x )= \det\{{\bf u}_{1}  (x ),\ldots,   {\bf u}_{N}  (x )\}.
 \label{eq:WROo}\end{equation}
\end{definition}

 Recall that     solutions  of equation \e{eq:H2} with asymptotics \e{eq:VV1} exist according  to Proposition~\ref{SRD}. Although they are not unique, according to Lemma~\ref{LI1}  the Wronskian ${\sf W}(x )$ does not depend (up to a numeration of eigenvalues $\z_{j}$) on a  specific choice of such solutions. In particular, we have
  \begin{equation}
 {\sf W} (x )= \det\{{\bf u}_{1}^{(-)}  (x; a_{-} ),\ldots,   {\bf u}_{n}^{(-)}  (x; a_{-} ),
  {\bf u}_{n+1}^{(+)}  (x; a_{+} ),\ldots,    {\bf u}_{N}^{(+)}  (x; a_{+} )\}
 \label{eq:WRO11}\end{equation}
 where $a_{+}$ ($a_{-}$) are sufficiently large positive (negative) numbers and the solutions ${\bf u}_{1}^{(\pm)}  (x; a_{\pm} )$ are constructed in Proposition~\ref{SRD}.

Of course, Definition~\ref{WRoo} applies if the matrices ${\bf L}_{0}(z)$ and ${\bf V} (x)$ are given by formulas \e{eq:J} and \e{eq:A}, respectively. In this case the  Wronskian ${\sf W}(x,z)$ depends analytically on the parameter $z\not\in \sigma (H_{0})$. Indeed, if additionally  $\Im z\neq 0$ for $N$ even and $\Re z\neq 0$ for $N$ odd, then this fact directly follows from the analyticity of the solutions $ {\bf u}_{j}^{(\pm)}  (x,z; a_{\pm} )$, $j=1,\ldots, N$
(see subs.~5.3). Moreover,  according to  Lemma~\ref{LI1}   the Wronskian ${\sf W}(x,z)$  is continuous (in contrast to the solutions  $ {\bf u}_{j}^{(\pm)}  (x,z; a_{\pm} )$) on the critical rays where $\Re \z_{l}=\Re \z_{j}$ for some
    $  \z_{l}\neq   \z_{j}$. Therefore its analyticity in a required region follows from Morera's theorem. Evidently, ${\sf W}(x,z)=0$  if and only if $z$ is an eigenvalue of the operator $H$.

  \medskip

 {\bf 6.2.}
  Let us return to  the trace formula \e{eq:tr} established so far
    for coefficients $v_k  (x)$, $k=1,\ldots, N$, with compact supports. Suppose that assumption \e{eq:V} holds. Then   condition  \e{eq:Vv}   is satisfied for both signs.  Let  us approximate $v_k  (x)$ by the cut-off functions $ \chi_{(-r, r)}(x) v_k  (x)$.
We denote by $u_{j}^{(-)}  (x, z; a_{-}, r )$, $j=1,\ldots, n$, and by $u_{j}^{(+)}  (x, z; a_+ , r )$, $j=n+1,\ldots, N$,   the solutions of equation  \e{eq:Heh} with the coefficients $\chi_{(-r, r) } v_k$  determined by Definition~\ref{DUa} and equality \e{eq:Qq}. Let us use formula \e{eq:WRO11} for    the   Wronskian
 ${\sf W}_{r}(x, z)$ for equation \e{eq:Heh} with  cut-off coefficients $ \chi_{(-r, r) }(x) v_k  (x)$.
    Then it follows from relation \e{eq:cut2} that
 \begin{equation}
\lim_{r\to\infty}{\sf W}_{r} (x,z)= {\sf W}  (x,z).
\label{eq:cutr}\end{equation}
In view of analyticity in $z$ of  these functions we also have
 \begin{equation}
\lim_{r\to\infty}\dot{\sf W}_{r} (x,z)= \dot{\sf W}  (x,z).
\label{eq:cutrd}\end{equation}

Set $H_{r}=H_{0}+V_{r}$ where the operator $V_{r}$ is defined by formula \e{eq:V1} with the coefficients   $ \chi_{(-r, r) }  v_k$.
Let us write down formula  \e{eq:resint2} for the operator $H_{r}$   and pass to the limit $r\to\infty$. By virtue of \e{eq:cutr} and \e{eq:cutrd}
  the right-hand side of  \e{eq:resint2} converges to the corresponding expression  for the operator $H $. It is possible to verify that the same is true for   the left-hand side of  \e{eq:resint2}. We shall not however dwell upon it and establish the trace formula in form
  \e{eq:tr}.

   Using the resolvent identity
   \[
   R(z)-R_{r}(z)=-\sum_{j=1}^N R_{r} (z) (V-V_{r})  R (z), \q R_{r}(z)= (H_{r}-z)^{-1},
   \]
   we see that for $z\not\in \sigma (H )$
 \[
\| R(z)-R_{r}(z)\|_{{\goth S}_1}  \leq C \| R_{0}(z) (V -V_{r} ) R_{0}(z)\|_{{\goth S}_1}
 \leq C_{1}\sum_{j=1}^N \| R_{0}(z) v_{j}(1-\chi_{r} )  \|_{{\goth S}_1}.
\]
 According to Proposition~\ref{TrCl} there is (for $N\geq 2$) the estimate
  \[
\| R_{0}(z) v_{j}(1-\chi_{r} )  \|_{{\goth S}_1}^2 \leq C
\int_{|x|\geq r}|v_{j}(x)|^2 (1+ x^2)^\alpha dx,\q \alpha>1/2,
\]
whence
\[
\lim_{r\to\infty} \| R_{r}(z) -R (z)\|_{{\goth S}_1}=0 .
\]
Thus, using trace formula  \e{eq:tr}  for  cut-off perturbations  $V_{r} $ and passing  to the limit $r\to\infty$, we deduce it for $V$. This leads to the following result.  Recall that  the normalized  Wronskian $\Delta(x,z)$ is defined by formula \e{eq:det2}.

  \begin{theorem}\label{MAIN}
 Under assumption \e{eq:V}  the trace formula  \e{eq:tr} holds  for all $z\not\in \sigma (H )$.
 \end{theorem}

   If  inclusion   \e{eq:TSS} holds,  then     equation  \e{eq:PD1} is satisfied for   a generalized perturbation determinant
   \begin{equation}
  \widetilde{D}_{z_{0}}(z)=
 \Det\big(I+ (z-z_{0})   R(z_{0}) V R_{0}(z )    \big)
  \label{eq:GC2}\end{equation}
  where $z_{0}\not\in \sigma (H)$.  It is easy to see that $\widetilde{D}_{z_{0}}(z)$
    is the usual perturbation determinant
   for the pair $R_{0}(z_{0}) $, $R (z_{0})$ at the point $(z-z_{0})^{-1}$. Of course, equation \e{eq:PD1}  for a function $\widetilde{D}(z)$ fixes  it up to a constant factor only.
   We note that for different ``reference points",  generalized perturbation determinants are connected by the formula
  $  \widetilde{D}_{z_{1}}(z)=   \widetilde{D}_{z_{0}}(z_{1})^{-1}\widetilde{D}_{z_0}(z)$. Moreover, if $ V R_{0}(z)\in {\goth S}_{1}$, then
  $  \widetilde{D}_{z_{0}}(z)=   D (z_{0})^{-1}  D(z)$  where   $D(z)$ is the  perturbation determinant  (see  formula
  \e{eq:PD}) for the pair $H_{0}$, $H$.

  Comparing equations \e{eq:tr} and \e{eq:PD1}  we see that for all $x\in \Bbb R$ and all   $z_{0}\not\in \sigma (H)$
  \[
  \Delta (x,z)= C(x_{0}, z_{0})\exp\big( -i^N \int_{x_{0}}^x v_{N}(y) dy\big) \widetilde{D}_{z_{0}}(z)
  \]
  where the constant $C(x_{0}, z_{0}) \neq 0 $ does not depend on $x$ and $z$.

 \medskip

  {\bf 6.3.}
  Suppose now that $v_{N}=0$.   Then the   Wronskian ${\sf W} (x, z)=:{\sf W} (z)$ does not depend on $x$,  and it is easy to deduce from Proposition~\ref{HE} that
   \begin{equation}
{\sf W}  ( z)= {\sf W}_{0}  ( z) \big(1+ O (  |z|^{-1/N})\big), \q | z| \to \infty.
  \label{eq:GD0}\end{equation}
  For the proof, it suffices to choose $a_{+}<0$, $a_{-}>0$ and use asymptotics \e{eq:Lx}
at $x=0$. As a side remark, we note that according to \e{eq:GD0} the set of complex eigenvalues of the operator $H$ is bounded.

  It follows from \e{eq:GD0} that the normalized Wronskian \e{eq:det2} satisfies the relation
     \begin{equation}
\Delta( z)= 1+ O (  |z|^{-1/N}), \q | z| \to \infty.
  \label{eq:GD}\end{equation}
 Since,   by  Proposition~\ref{TrCl}, $V R_{0}(z)\in{\goth S}_{1}$,   the  perturbation determinant  is correctly defined by formula
  \e{eq:PD} and (see book \cite{GoKr})
  \begin{equation}
\lim_{| \Im z | \to\infty} \Det\big(I+V R_{0}(z)\big) =1.
\label{eq:PD2}\end{equation}
  Comparing equations  \e{eq:tr} and \e{eq:PD1}, we obtain  that
\begin{equation}
\Delta(z)=C \Det\big(I+V R_{0}(z)\big)
\label{eq:PDC}\end{equation}
for some constant $C$.
 Moreover, taking into account relations    \e{eq:GD} and \e{eq:PD2}, we see that $C=1$ in   \e{eq:PDC}.

 Let us formulate      the result obtained.

  \begin{theorem}\label{MAIN1}
  Suppose that  $v_{N}=0$ and that the coefficients $v_{j}$, $j= 1,\ldots, N-1$, satisfy  assumption   \e{eq:V}.   Then $ \Delta(x, z)= :\Delta(z)$ does not depend on $x$ and  relation  \e{eq:PDd} is true  for all $z\not\in \sigma (H )$.
\end{theorem}

    \medskip

    {\bf 6.4.}
    Finally, we note that for a derivation of the trace formula  \e{eq:tr} the approximation of $v_{j}$ by cut-off functions $\chi_{(-r,r)}v_{j}$ is not really necessary. From the very beginning, we could work with functions $v_{j}$  satisfying assumption \e{eq:V} only. Then the text of Sections~2, 3 and 4 remains unchanged if, for all $j=1,\ldots, N$, the functions $e^{\z_{j}x}{\bf p}_{j}$ are replaced for $x<\!\!<0$ by ${\bf u}_{j}^{(-)}(x ;a_{-})$ where $a_{-}$ is a sufficiently big negative number and they are replaced for $x>\!\!>0$ by ${\bf u}_{j}^{(+)}(x ; a_{+})$ where $a_{+}$ is a sufficiently big positive number. In particular, the definition of the transition matrices in subs.~2.4 can be given in terms of the solutions ${\bf u}_{j}^{(\pm)}(x ; a_{\pm})$.

    However, a preliminary consideration of coefficients $v_{j}$ with compact supports seems to be intuitively more clear.

    \end{document}